\documentclass[10pt,a4paper]{article}
\usepackage{authblk}
\usepackage{amsthm}
\usepackage{amsmath}
\usepackage{amssymb}
\usepackage{tikz-network}
\usepackage{relsize}
\usepackage{multirow}
\usepackage{xcolor}
\usepackage[verbose]{hyperref}
\hypersetup{colorlinks=false,allbordercolors=blue,pdfborderstyle={/S/U/W 1}}
\usepackage{enumitem}

\theoremstyle{plain}% Theorem-like structures provided by amsthm.sty
\newtheorem{theorem}{Theorem}[section]
\newtheorem{lemma}[theorem]{Lemma}
\newtheorem{corollary}[theorem]{Corollary}
\newtheorem{proposition}[theorem]{Proposition}

\newtheorem{remark}[theorem]{Remark}
\theoremstyle{definition}
\newtheorem{definition}[theorem]{Definition}
\newtheorem{example}[theorem]{Example}

\DeclareMathOperator{\CS}{CS}
\DeclareMathOperator{\SC}{SC}
\newcommand{\vvee}{% bitwise logical and
	\hbox{%
		$\lor$%
		\hspace{-.47em}%
		\raisebox{0.60em}{%
			\rotatebox{-66}{%
				\vrule width .30em height .45pt depth 0pt%
			}%
			\hspace{.065em}%
			\rotatebox{66}{%
				\hspace{-.29em}%
				\vrule width .30em height .45pt depth 0pt%
			}%
		}%
		\hspace{.20em}%
	}%
}

\DeclareMathOperator*{\hjoin}{\vvee}
\begin{document}

\title{On the structural aspects of Convex Subhyperlattices under different orders}

\author[1]{ N Aishwarya Nayak} 
\author[2]{Bhargava K}
 \author[3]{Pallavi Panjarike}
 \author[1]{Kuncham Syam Prasad}
 \author[3]{Tapatee Sahoo}
 \author[1]{Harikrishnan Panackal}

\affil[1]{Manipal Institute of Technology,\\ Manipal Academy of Higher Education, Manipal, India\\}
\affil[2]{Department of Mathematics,Mangalore University\\Mangalagangothri, Karnataka, India\\}
\affil[3]{Manipal Institute of Technology Bengaluru,\\ Manipal Academy of Higher Education, Manipal, India\\}
	
\maketitle

\begin{abstract}
Nayak et al. (2025) introduced the concept of convex subhyperlattices of a strong join hyperlattice $\mathfrak{L}$ and demonstrated that the collection of all convex subhyperlattices, $\CS(\mathfrak{L})$, together with empty set, forms a lattice. In particular, for $P$-hyperlattices and Nakano hyperlattices, they obtained a unique representation of convex subhyperlattices in terms of hyperideals and hyperfilters. In the present work, we extend these studies by proving that $\CS(\mathfrak{L}) \cup \{\emptyset\}$, under set inclusion, forms an algebraic lattice. Moreover, we introduce two new orders on $\CS(\mathfrak{L})$ and examine their properties, including conditions for equivalence of these orders. In the class of $P$-hyperlattices, we characterize convex subhyperlattices in terms of initial and final segments, providing a structural insight into their formation. When convexity is extended to arbitrary join hyperlattices, many results known for $P$-hyperlattices and Nakano hyperlattices fail to hold. To address these limitations of convexity in arbitrary hyperlattices, we introduce the notion of strongly convex sets, denoted $\SC(\mathfrak{L})$, which naturally generalizes the concept while preserving its fundamental properties and extends the applicability of results to a wider setting. We establish several results for $\SC(\mathfrak{L})$ that do not hold for $\CS(\mathfrak{L})$.
\end{abstract}

\subsection*{Keywords}
hyperlattice, poset, convex set, cover, lattice homomorphism.

\subsection*{ Mathematics Subject Classification:}
 06B10, 06B23, 06B75

\section{Introduction and Preliminaries}\label{sec1}
Hyperoperations serve as a broad generalization of binary operations. Building on this idea, Konstantinidou and Mittas \cite{1977} extended the classical notion of lattices by replacing the binary join and meet operations with corresponding hyperoperations. In \cite{Pallavi2022}, the authors introduced and investigated various generalizations of prime hyperideals, such as 2-absorbing and primary 2-absorbing hyperideals, within the framework of meet hyperlattices.  Bideshki et al. \cite{Bideshki2018} introduced the ideas of strong meet hyperlattices, hyperideals, and hyperfilters in strong meet hyperlattices. They also studied distributivity and dual distributivity in meet hyperlattices. Ameri et al. \cite{Ameri2014} examined residual mappings in meet hyperlattices through the framework of principal downsets and established specific conditions under which a mapping between two hyperlattices corresponds to a residual mapping. They further explored the structure and properties of convex subsets in meet hyperlattices. Later, Nayak et al. \cite{Nayak2025} introduced two kinds of hyperfilters in join hyperlattices and explored their relationships within $P$-hyperlattices and Nakano hyperlattices. They also examined convex subhyperlattices and provided a unique representation for them in terms of hyperideals and hyperfilters within $P$-hyperlattices and Nakano hyperlattices. Lavanya and Parameshwara Bhatta~\cite{lavanya1996new} investigated the collection of all convex sublattices of a lattice $L$ under a partial order distinct from inclusion. They showed that the lattice of convex sublattices under this new order and the lattice under consideration belong to the same equational class. Extending their ideas, Bhargava and Hegde~\cite{bhargava2025} introduced the concept of strongly convex subtrellises in trellises, a generalization of lattices, and generalized several corresponding results. Since hyperlattices are also natural generalizations of lattices, it is meaningful to ask whether these developments can be extended to the framework of hyperlattices. Nayak et al. \cite{Nayak2025} established that the set $\CS(\mathfrak{L})$ of all convex subhyperlattices of a join hyperlattice $\mathfrak{L}$, including the empty set, forms a lattice under inclusion. In the present study, we demonstrate that $\CS(\mathfrak{L}) \cup \{\emptyset\}$ is an algebraic lattice. Furthermore, we introduce a new order on $\CS(\mathfrak{L})$, defined for $X,Y \in \CS(\mathfrak{L})$ such that $X \leq_b Y$ if and only if, for every element $x \in X$, there exists $y \in Y$ such that $x \preceq y$, and for every $y \in Y$, there exists $x \in X$ such that $x \preceq y$. Under this relation, $\CS(\mathfrak{L})$ forms a partially ordered set. We also define another order on $\CS(\mathfrak{L})$ using hyperideals and hyperfilters, and show that these two orders are equivalent in $P$-hyperlattices, though not necessarily in general, as illustrated through examples. Although convexity is well defined in strong hyperlattices, its direct extension to arbitrary hyperlattices is inadequate for preserving key results. The absence of unique representation of convex subhyperlattices motivates the need for the introduction of strongly convex sets. In the last section, we introduce strongly convex sets and discuss some of their properties.
\newline
\newline

Throughout this paper we use the following notations:
\begin{itemize}
	\item $\vee$ is the join of two elements in a lattice $L$.
	\item $\hjoin$ is the hyperoperation join.
	\item $\bigvee$ is the classical join of a subset of a lattice $L$.
	\item $\preceq$ is the partial order on the hyperlattice $\mathfrak{L}$.
	\item $\prec$ is the covering relation in a classical lattice.
\end{itemize}
\begin{definition}[\cite{Davvaz2007}]
	Let $ H $ be a non-empty set and let $\mathcal{P}^*(H)$ be the set of all non-empty subsets of $ H $.  A hyperoperation on $ H $ is a map $\circ: H \times H\rightarrow \mathcal{P}^*(H) $ and $(H,\circ)$ is called a hypergroupoid.
\end{definition}  
\noindent If $ A $ and $ B $ are non-empty subsets of $ H $ and $\hjoin$ is a hyperoperation, then
$$A\hjoin B =\bigcup_{a\in A,b\in B} a\hjoin b,~~~A \hjoin  x=A \hjoin  \{x\}~~and~~x\hjoin B=\{x\}\hjoin B.$$
\begin{definition}[\cite{SoltaniLashkenari2015}]
	Let $L$ be a non-empty set, $\hjoin:L\times L \rightarrow \mathcal{P}^*(L)$ be a hyperoperation, and $\wedge: L\times L \rightarrow L$ be a binary operation. Then $(L,\hjoin,\wedge)$ is a join hyperlattice if for all $x,y,z\in L$ the following conditions hold: 	
	\begin{enumerate}[label=\roman*]
		\item $x\in x\hjoin x$ and $x=x\wedge x$.
		\item $x\hjoin(y\hjoin z)=(x\hjoin y)\hjoin z$ and $x\wedge(y\wedge z)=(x\wedge y)\wedge z$.
		\item $x\hjoin y=y\hjoin x$ and $ x\wedge y=y\wedge x$.
		\item  $x\in x\wedge(x\hjoin y)\cap x\hjoin(x\wedge y )$.
	\end{enumerate}
	We say $L$ is a strong join hyperlattice if whenever $y\in x\hjoin y$ then $x=x\wedge y $.
\end{definition}

\noindent Throughout, $\mathfrak{L}$ denotes a strong join hyperlattice unless otherwise specified, and $\hjoin$ and $\vee$ denote the join hyperoperation and the classical join operation respectively.
\begin{definition}[\cite{Rasouli2010a}] \label{defn: partialorderonHL}
	Define the relation $ \preceq $ on  $ \mathfrak{L} $ as follows: $$  l_1\preceq l_2~~  \text{if and only if}~~   l_1\wedge l_2 =l_1 ~\text{if and only if}~ l_2 \in l_1 \hjoin l_2.$$ Then  $( \mathfrak{L},\preceq) $ is poset.
\end{definition}
Henceforth, $\preceq$ will denote a partial order on the strong hyperlattice $\mathfrak{L}$.

A non-empty subset $H$ of $\mathfrak{L}$ is a subhyperlattice if $H$ is closed under the hyperoperation $\hjoin$ and the operation $\wedge$ in $\mathfrak{L}$~\cite{Rasouli2010}. A non-empty subset $J $ of $ \mathfrak{L} $ which is closed under $\hjoin$ is a hyperideal if $a \in \mathfrak{L}, b \in J$ with $a \preceq b$ implies that $a \in J$~\cite{SoltaniLashkenari2015}. A non-empty subset $F$ of $\mathfrak{L}$ which is closed under $\wedge$ is a hyperfilter of 
\begin{enumerate}[label=\roman*]
	\item \textit{type-I}, if whenever $a \in F, x \in \mathfrak{L}$ with $a \preceq x$ implies $x \in F$~\cite{SoltaniLashkenari2015}.
	\item \textit{type-II} if $a_1 \in F, x \in \mathfrak{L}$ with $a_1 \preceq x$ implies $a_1 \hjoin x \subseteq F$~\cite{Nayak2025}.
	\item 	\textit{type-III} if $a \in F, x \in \mathfrak{L}$ implies that $a \hjoin x \subseteq F$~\cite{Nayak2025}.
\end{enumerate}

\begin{definition}[\cite{Rasouli2010a}]
	Let $A \subseteq \mathfrak{L}$. The intersection of all hyperideals of $\mathfrak{L}$ containing~$A$ is called the hyperideal generated by $A$ and it is denoted by $\langle A]$. If $A=\{a\}$, then we call the hyperideal generated by $A$ a principal hyperideal and denote it by $\langle a]$.
\end{definition}
\begin{definition}
	The intersection of all type-III hyperfilters of $\mathfrak{L}$ containing~$A$ is called the type-III hyperfilter generated by $A$ (or simply, the hyperfilter generated by $A$). 
\end{definition}
In the sequel, we use the notation $\langle A]$ to denote the hyperideal generated by $A$ and $[A\rangle$ to denote the (type-III) hyperfilter generated by $A$. \newline
The set of all hyperideals of $\mathfrak{L}$ forms a lattice $I(\mathfrak{L})$ under inclusion, where the meet and join operations are given as: for $J_1, J_2 \in I(\mathfrak{L})$, $J_1 \wedge J_2 = J_1 \cap J_2$ and $J_1 \vee J_2 = \langle J_1 \cup J_2]$. Similarly, the set of all type-III hyperfilters forms a lattice $F(\mathfrak{L})$ under reverse inclusion where the meet and join operations are given as: for $J_1, J_2 \in F(\mathfrak{L})$, $J_1 \wedge J_2 = [J_1 \cup J_2 \rangle $ and $J_1 \vee J_2 =  J_1 \cap J_2$.
\begin{theorem} [\cite{Konstantinidou1993}]
	Let $(L,\wedge,\vee)$ be lattice and $P$ be a non-empty subset of $ L.$
	We define a hyperoperation $\hjoin\limits^{P}$ by 
	$$l_{1} \hjoin\limits^P l_{2}=l_{1}\vee l_2 \vee P=\{l_{1}\vee l_{2}\vee p~|~p\in P\},$$
	
	\noindent for all $l_1, l_2 \in L$.	Then $(L,\hjoin\limits^P,\wedge)$	is a join hyperlattice if and only if for each $l_1 \in L$ there exists a $p\in P $ such that $p\preceq l_{1}$. In this case we call $L$ a $P$-hyperlattice. Also, $L$ is a strong hyperlattice.
\end{theorem}  

\begin{definition}[\cite{Rasouli2010}]
	Let $(L, \wedge , \vee)$ be a modular lattice. For any $a,b \in L$, define the Nakano hyperoperation $a \sqcup b=\{\, z \in L : z \vee a=z \vee b = a \vee b \,\}$. Then $(L, \sqcup, \wedge)$ is a strong join hyperlattice.
	
\end{definition}
We call the hyperlattice $(L, \sqcup, \wedge)$ a Nakano join hyperlattice.

\begin{lemma}[\cite{Nayak2025}]
	Intersection of hyperfilters of one type is again a  hyperfilter of the same type, provided it is non-empty.
\end{lemma}

\begin{proposition}[\cite{Nayak2025}]~
	\begin{enumerate}[label=\roman*]
		\item 	In a join $P$-hyperlattice, all the three types of hyperfilters are equivalent.
		\item 	In a Nakano join hyperlattice, the type-II and type-III hyperfilters are equivalent.
	\end{enumerate}
\end{proposition}

\begin{definition}[\cite{Nayak2025}]
	A subhyperlattice $K$ of $\mathfrak{L}$ is a convex subhyperlattice if for any $p_1,p_2 \in K,~ p_3 \in \mathfrak{L}$ with $p_1 \preceq p_3 \preceq p_2$ implies that $p_3 \in K$.
	
\end{definition}

\begin{proposition}[\cite{Nayak2025}] \label{proposition5}
	For any $\mathfrak{L}$,
	\begin{enumerate}[label=\roman*]
		\item All three types of hyperfilters of $\mathfrak{L}$ are convex subsets of $\mathfrak{L}$.
		\item 	Every hyperideal of $\mathfrak{L}$ is a convex subhyperlattice.
		\item 	Every type-III hyperfilter  of  $\mathfrak{L}$ is a convex subhyperlattice.
	\end{enumerate}	
\end{proposition}

\begin{proposition}[\cite{Nayak2025}]  \label{proposition3}
	Let $\{C_i\}_{i \in \tau}$ be a collection of convex subhyperlattices of $\mathfrak{L}$. Then $C_o = \bigcap \limits_{i \in \tau} C_i$ is a convex subhyperlattice provided $C_o$ is non-empty.
\end{proposition}

\begin{lemma}[\cite{Nayak2025}]
	Let $CS(\mathfrak{L})$ be the set of all convex subhyperlattices of~  $\mathfrak{L}$. Then $(CS(\mathfrak{L}) \cup \{\emptyset\}, \subseteq)$ forms a lattice under inclusion order.
\end{lemma}

\begin{theorem}[\cite{Nayak2025}]\label{thmC}
	Let $I$ be a hyperideal and $F$ be a type-III hyperfilter of $\mathfrak{L}$. Then $I \cap F$ is a convex subhyperlattice provided it is non-empty.
\end{theorem}

\begin{theorem}[\cite{Nayak2025}] \label{thmP}
	In a join $P$-hyperlattice, any convex subhyperlattice can be  uniquely written as the intersection of a hyperideal and a hyperfilter.
\end{theorem}

\section{On the lattice $(\CS(\mathfrak{L}) \cup \{\emptyset\},\subseteq)$}

The authors in \cite{Nayak2025}, have established that $(\CS(\mathfrak{L}) \cup \{\emptyset\}, \subseteq)$ is a lattice. In this section, we explore various properties of the lattice $(\CS(\mathfrak{L}) \cup \{\emptyset\}, \subseteq)$ in detail. We derive the expression for the convex subhyperlattice generated by a non-empty subset of a join hyperlattice, and characterize the compact elements of $(\CS(\mathfrak{L}) \cup \{\emptyset\}, \subseteq)$. Using these results, we show that the lattice $(\CS(\mathfrak{L}) \cup \{\emptyset\}, \subseteq)$ is algebraic.

\begin{definition}
	For a non-empty subset $A$ of a join hyperlattice $\mathfrak{L}$, the convex subhyperlattice generated by $A$ is denoted by $\langle A \rangle$ and is defined as the intersection of all convex subhyperlattices containing $A$.
\end{definition}

For a non-empty subset $A$ of a  $P$-hyperlattice $\mathfrak{L}$, the convex subhyperlattice generated by $A$ is given by  $\langle A \rangle = \langle A] \cap [A\rangle$, where $\langle A]$ denotes the hyperideal generated by $A$ and $[A \rangle $ denotes the hyperfilter generated by $A$. \newline If $A=\{a\}$, we write $\langle a \rangle$ instead of $\langle \{a\} \rangle$.

\begin{remark}For non-empty subsets $H_1$ and $H_2$ of a  $P$-hyperlattice $\mathfrak{L}$,
	\begin{enumerate}[label=\roman*]
		\item We observe that  $ \langle H_1] \cap [H_1\rangle$ is the smallest convex subhyperlattice containing $H_1$. Indeed, $ \langle H_1] \cap [H_1\rangle$ is itself a convex subhyperlattice containing $H_1$. If $G_1$ is any convex subhyperlattice containing $H_1$, then $\langle H_1] \subseteq\langle G_1]$ and $[H_1\rangle \subseteq [G_1\rangle$. Hence, $ \langle H_1] \cap [H_1\rangle \subseteq \langle G_1] \cap [G_1\rangle=G_1$.
		\item $H_1$ is a convex subhyperlattice if and only if $\langle H_1 \rangle=H_1$.
		\item If $H_1 \subseteq H_2$, then $\langle H_1 \rangle \subseteq \langle H_2 \rangle$.
	\end{enumerate}
\end{remark}
\begin{remark}
	Let $(\mathfrak{L},\hjoin,\wedge)$ be any join hyperlattice and $a \in \mathfrak{L}$. Then $\{a\}$ is a convex subhyperlattice if and only if $a \hjoin a =\{a\}$.
\end{remark}

\begin{proposition}
	For any $\mathfrak{L}$, $(\CS(\mathfrak{L}) \cup \{\emptyset\},\subseteq)$ forms a complete lattice where for $X,Y \in \CS(\mathfrak{L}),$
	\begin{align*}
		X \wedge Y &= X \cap Y \\ 
		X \vee Y &= \langle X \cup Y \rangle
	\end{align*}
\end{proposition}
\begin{proof}
	This follows from the fact that the intersection of any collection of convex subhyperlattices of $\mathfrak{L}$ is either empty or a convex subhyperlattice.
\end{proof}

\begin{lemma}\label{lem:unionconvex}
	Let $(\mathfrak{L}, \hjoin, \wedge)$ be a strong join hyperlattice. If $\{X_\alpha\}$ is a chain in $(\CS(\mathfrak{L}) \cup \{\emptyset\}, \subseteq)$, then $X= \bigcup\limits_{\alpha} X_\alpha$ is a convex subhyperlattice of $\mathfrak{L}$.
\end{lemma}
\begin{proof}
	Let $x, y \in X$  and $a \in \mathfrak{L}$ such that $x \preceq a \preceq y$. Since $\{X_\alpha\}$ is a chain, we get $x,y \in X_\beta$ for some $\beta$. By the convexity of $X_\beta$ we then have $a \in X_\beta \subseteq X$. Hence $X$ is a convex set. Next, let $a, b \in X$. Then $a, b \in X_\alpha$ for some $\alpha$. Now, $X_\alpha$ is a subhyperlattice implies that $a \wedge b \in X_\alpha \subseteq X$ and $a \hjoin b \subseteq X_\alpha  \subseteq X$. Therefore $X$ is a subhyperlattice.
\end{proof}

\begin{proposition}\label{prop:convexgeneration}
	For a non-empty subset $A$ of $\mathfrak{L}$,
	\begin{equation*}
		\langle A \rangle  = \bigcup \{\langle F \rangle  \mid F \subseteq A, F~ \text{is finite}\}.	
	\end{equation*}
\end{proposition}
\begin{proof}
	The proof is by transfinite induction on $|A|$. If $A$ is finite, then the result holds trivially. Now, suppose that $A$ is infinite and the result holds for all sets of cardinality smaller than $|A|$. We well-order the set $A$ in the order type of its cardinal. Then for each $a \in A$, the initial segment $S(a) =\{x \in A \mid x<a\}$ determined by $a$ has cardinality strictly less than $|A|$. Furthermore, as $A$ is infinite, its cardinal is a limit ordinal and hence it follows that
	\begin{equation}\label{eq:1}
		A= \bigcup\limits_{a \in A} S(a)
	\end{equation}
	Let $\mathcal{C} =\{\langle S(a) \rangle \mid a \in A\}$. Let $a,b \in A$. Since $A$ is a chain with respect to $\leq$, we get $a \leq b$ or $b \leq a$. Without loss of generality, assume that $a \leq b$. Then $S(a) \subseteq S(b)$ and hence, $\langle S(a) \rangle \subseteq \langle S(b) \rangle$. This shows that $\mathcal{C}$ is a chain with respect to $\subseteq$. Let $\textbf{C} =\bigcup \mathcal{C} = \bigcup\limits_{a\in A} \langle S(a) \rangle$. Then by Lemma~\ref{lem:unionconvex}, $\textbf{C}$ is a convex subhyperlattice. Now, for every $a \in A$, since $S(a) \subseteq \langle S(a) \rangle \subseteq \textbf{C}$, it follows from equation (\ref{eq:1}) that $A \subseteq \textbf{C}$. Hence $\langle A \rangle \subseteq \langle \textbf{C} \rangle=\textbf{C}$ because $\textbf{C} \in \CS(\mathfrak{L})$. Also for any $a \in A$, since $S(a) \subseteq A$, we get $\langle S(a) \rangle \subseteq \langle A \rangle$ which implies $\textbf{C} \subseteq \langle A \rangle$. Thus we get $\langle A \rangle = \textbf{C}$.  By induction hypothesis, for each $a \in A$, we have 
	\begin{equation*}
		\langle S(a) \rangle = \bigcup \{ \langle F \rangle \mid F \subseteq S(a), F~\text{is finite} \}
	\end{equation*}
	Thus, $\langle A \rangle = \bigcup \limits_{a \in A} \langle S(a) \rangle = \bigcup  \{\langle F \rangle  \mid F \subseteq A, F~ \text{is finite}\}$, as any finite subset $F$ of $A$ is contained in $S(a)$ for some $a \in A$.
\end{proof}

\begin{definition}[\cite{Gratzer2011}]
	An element $p$ of a lattice $L$ is compact if whenever $p \leq \bigvee S$ for a subset $S$ of $L$, then $p \leq \bigvee F$ for a finite subset $F$ of $S$. 
\end{definition}
\begin{definition}[\cite{Gratzer}]
	A complete lattice $L$ is said to be algebraic (compactly generated) if every element of $L$ is a join of some compact elements of $L$. 
\end{definition}

\begin{lemma}\label{lem3}
	Let $S=\{x_1,x_2,\dots, x_n\}$ be a non-empty finite subset of $\mathfrak{L}$. Then $\langle S \rangle= \mathlarger{\bigvee}\limits_{i=1}^{n}\langle x_i \rangle$.
\end{lemma}
\begin{proof}
	We know that $x_i \in \langle x_i \rangle$ for each $i,~ 1 \leq i \leq n$. This gives $S \subseteq \mathlarger{\bigcup}\limits_{i=1}^{n}\langle x_i \rangle$ and hence $\langle S \rangle \subseteq \Biggl \langle\mathlarger{\bigcup}\limits_{i=1}^{n}\langle x_i \rangle \Biggr \rangle = \mathlarger{\bigvee}\limits_{i=1}^{n}\langle x_i \rangle$. Since $x_i \in S$, we have $\langle x_i \rangle \subseteq \langle S \rangle$ for each $i$. This implies $\mathlarger{\bigcup}\limits_{i=1}^{n}\langle x_i \rangle \subseteq \langle S \rangle$. Thus $\mathlarger{\bigvee}\limits_{i=1}^{n}\langle x_i \rangle= \Biggl \langle\mathlarger{\bigcup}\limits_{i=1}^{n}\langle x_i \rangle \Biggr \rangle \subseteq\langle S \rangle$. Hence $\langle S \rangle= \mathlarger{\bigvee}\limits_{i=1}^{n}\langle x_i \rangle$.
\end{proof}

\begin{proposition}
	For any $(\mathfrak{L}, \hjoin, \wedge)$, the complete lattice $(\CS(\mathfrak{L})\cup\{\emptyset\}, \subseteq)$ is algebraic. The set of all compact elements of this lattice is 
	\begin{equation*}
		\{\langle F \rangle \mid F ~\text{is a finite subset of}~ \mathfrak{L}\}
	\end{equation*}
\end{proposition}
\begin{proof}
	Let $F=\{f_1,f_2,\cdots, f_n\}$ be a finite subset of $\mathfrak{L}$, for some $n \in \mathbb{N}$. We show that $\langle F \rangle$ is compact in $(\CS(\mathfrak{L})\cup\{\emptyset\}, \subseteq)$. Suppose $\langle F \rangle \subseteq \mathlarger{\bigvee}\limits_{\alpha \in \mathcal{I}} X_\alpha = \Biggl \langle \mathlarger{\bigcup}\limits_{\alpha \in \mathcal{I}} X_\alpha \Biggr \rangle$ for some sub collection $\{X_\alpha \mid \alpha \in \mathcal{I}\}$ in $\CS(\mathfrak{L})$. By Proposition~\ref{prop:convexgeneration}, we have 
	\begin{equation*}
		\Biggl \langle \mathlarger{\bigcup}\limits_{\alpha \in \mathcal{I}} X_\alpha \Biggr \rangle = \mathlarger{\bigcup} \{\langle S \rangle \mid S ~\text{is finite subset of}~ \mathlarger{\bigcup}\limits_{\alpha \in \mathcal{I}} X_\alpha\}
	\end{equation*}
	Since $f_i \in F \subseteq \langle F \rangle$ for each $i$, there exists a finite subset $S_i$ of $\mathlarger{\bigcup}\limits_{\alpha \in \mathcal{I}} X_\alpha$ such that $f_i \in \langle S_i \rangle$. Each $S_i$ is a finite subset of $\mathlarger{\bigcup}\limits_{\alpha \in \mathcal{I}} X_\alpha$. So, there exists a finite subset $\mathcal{I}_i \subseteq \mathcal{I}$ such that $S_i \subseteq \mathlarger{\bigcup}\limits_{\alpha \in \mathcal{I}_i} X_\alpha$. Put $\Lambda_F = \mathlarger{\bigcup}\limits_{j=1}^{n} \mathcal{I}_j$. Then $\Lambda_F$ is a finite subset of $\mathcal{I}$. We have $f_j \in \langle S_j \rangle$, which implies
	\begin{eqnarray*}
		F \subseteq \mathlarger{\bigcup}\limits_{j=1}^{n}\langle S_j \rangle &\subseteq& \Biggl \langle \mathlarger{\bigcup}\limits_{j=1}^{n}S_j \Biggr \rangle \\ &\subseteq& \Biggl \langle \mathlarger{\bigcup}\limits_{j=1}^{n} \mathlarger{\bigcup}\limits_{\alpha \in \mathcal{I}_j}X_\alpha \Biggr \rangle \\ &=& \Biggl \langle \mathlarger{\bigcup}\limits_{\alpha \in \Lambda_F}X_\alpha \Biggr \rangle. 
	\end{eqnarray*}
	Hence, $\langle F \rangle \subseteq \Biggl \langle \mathlarger{\bigcup}\limits_{\alpha \in \Lambda_F}X_\alpha \Biggr \rangle= \mathlarger{\bigvee}\limits_{\alpha \in \Lambda_F}X_\alpha$. Therefore, we get $\langle F \rangle$ is a compact element in $(\CS(\mathfrak{L})\cup\{\emptyset\}, \subseteq)$. Now, let $X \in \CS(\mathfrak{L})$ be any compact element. Clearly, $X \subseteq  \mathlarger{\bigcup}\limits_{x\in X} \langle x \rangle \subseteq \Biggl \langle \mathlarger{\bigcup}\limits_{x\in X}\langle x \rangle \Biggr\rangle = \mathlarger{\bigvee}\limits_{x\in X}\langle x \rangle$. Since $X$ is compact, there exist finitely many elements $x_1, x_2, \dots, x_n \in X$ such that 
	\begin{equation}\label{eq:2}
		X \subseteq \mathlarger{\bigvee}\limits_{i=1}^{n} \langle x_i \rangle.
	\end{equation} Let $F=\{x_1,x_2,\dots,x_n\}$. By Lemma~\ref{lem3}, $\langle F \rangle = \mathlarger{\bigvee}\limits_{i=1}^{n}\langle x_i \rangle$. From equation (\ref{eq:2}), we get $X \subseteq \langle F \rangle$. Since $F\subseteq X$, we have $\langle F \rangle \subseteq \langle X \rangle =X$. Hence $X=\langle F \rangle$ where $F$ is a finite subset of $\mathfrak{L}$. Thus the set of all compact elements of the complete lattice $(\CS(\mathfrak{L})\cup\{\emptyset\}, \subseteq)$ is $\{\langle F \rangle \mid F ~\text{finite subset of}~ \mathfrak{L}\}$
	
	Finally, we show that $(\CS(\mathfrak{L})\cup\{\emptyset\}, \subseteq)$ is algebraic. Let $X \in \CS(\mathfrak{L})$. By Proposition~\ref{prop:convexgeneration}, we have 
	\begin{equation*}
		X =\mathlarger{\bigcup} \{\langle F \rangle  \mid F \subseteq X, F~ \text{is finite}\}
	\end{equation*}
	Then, 
	\begin{eqnarray*}
		\mathlarger{\bigvee} \{\langle F \rangle  \mid F \subseteq X, F~ \text{is finite}\} &=& \Biggl \langle \mathlarger{\bigcup} \{\langle F \rangle  \mid F \subseteq X, F~ \text{is finite}\} \Biggr \rangle\\ &=& \langle X \rangle = X ~(\text{because}~ X~ \text{is convex subhyperlattice})
	\end{eqnarray*}
	This shows that $X$ is a join of some compact elements in $(\CS(\mathfrak{L})\cup\{\emptyset\}, \subseteq)$.
\end{proof}

\begin{example}\label{eg1}
	Consider the lattice $\mathfrak{L}=\{0,a,b,c,d,e,f,h,1\}$ shown in Figure~\ref{fig1}. The $P$-hyperlattice obtained by taking $P=\{0,c\}$ is given in the Figure~\ref{fig2}. In the diagrammatic representation, the meet of two elements is the lattice element occurring in the same position as their meet in the lattice. The hyperoperation join of two elements is the set occurring in the same position as their join in the lattice. For instance, $0 \hjoin 0= \{0,c\}, a \hjoin e =\{1\}, c \hjoin f =\{f\}$.
	
	\begin{figure}[th]
		\begin{minipage}[b]{0.45\linewidth}
			\begin{center}
				
				\begin{tikzpicture}[scale=0.8]
					\draw[fill=black] (0,0) circle(2pt);
					\draw[fill=black] (-3,1.5) circle(2pt);
					\draw[fill=black] (-1.5,1.5) circle(2pt);
					\draw[fill=black] (0,1.5) circle(2pt);
					\draw[fill=black] (1.5,1.5) circle(2pt);
					\draw[fill=black] (3,1.5) circle(2pt);
					\draw[fill=black] (-1.5,3) circle(2pt);
					\draw[fill=black] (1.5,3) circle(2pt);
					\draw[fill=black] (0,4.5) circle(2pt);
					
					\node at (0,-0.3){$0$};
					\node at (-3.3,1.5){$a$};
					\node at (-1.8,1.5){$b$};
					\node at (-0.3,1.5){$c$};
					\node at (1.8,1.5){$d$};
					\node at (3.3,1.5){$e$};
					\node at (-1.8,3){$f$};
					\node at (1.8,3){$h$};
					\node at (0,4.8){$1$};
					
					\draw[thick]  (0,0)--(-3,1.5)--(-1.5,3)--(-1.5,1.5)--(0,0)--(0,1.5)--(-1.5,3)--(0,4.5)--(1.5,3)--(0,1.5);
					\draw[thick] (1.5,3)--(1.5,1.5)--(0,0)--(3,1.5)--(1.5,3);
					
				\end{tikzpicture}
				
				\caption{Lattice $\mathfrak{L}$ from Example~\ref{eg1} }\label{fig1}
			\end{center}
		\end{minipage}
		\begin{minipage}[b]{0.45\linewidth}
			\begin{center}
				
				\begin{tikzpicture}[scale=0.7]
					\node (0) at (0,-1) {\begin{tabular}{c}
							$ 	0 $\\$ \{0,c\} $
					\end{tabular}};

					\node (a) at (-3,2){\begin{tabular}{c}
							$ a $\\$ \{a,f\} $ 
					\end{tabular}};
					\node (b) at (-1.5,2){\begin{tabular}{c}
							$ b$\\$ \{b,f\} $ 
					\end{tabular}};
					\node (c) at (0,2){\begin{tabular}{c}
							$ c$\\$ \{c\} $ 
					\end{tabular}};
					\node (d) at (1.5,2){\begin{tabular}{c}
							$ d $\\$ \{d,h\} $ 
					\end{tabular}};
					\node (e) at (3,2){\begin{tabular}{c}
							$ e $\\$ \{e,h\} $ 
					\end{tabular}};
					\node (f) at (-1.5,5){\begin{tabular}{c}
							$  f$\\$ \{f\} $ 
					\end{tabular}};
					\node (h) at (1.5,5){\begin{tabular}{c}
							$h$\\$ \{h\} $ 
					\end{tabular}};
					\node (1) at (0,8){\begin{tabular}{c}
							$ 1 $\\$ \{1\} $ 
					\end{tabular}};
					
					\draw[thick]  (0)--(a)--(f)--(b)--(0)--(c)--(f)--(1)--(h)--(c);
					\draw[thick] (h)--(d)--(0)--(e)--(h);
					
				\end{tikzpicture}
				\caption{ Hyperlattice $\mathfrak{L}$ from Example~\ref{eg1} }\label{fig2}
			\end{center}
		\end{minipage}
	\end{figure}

	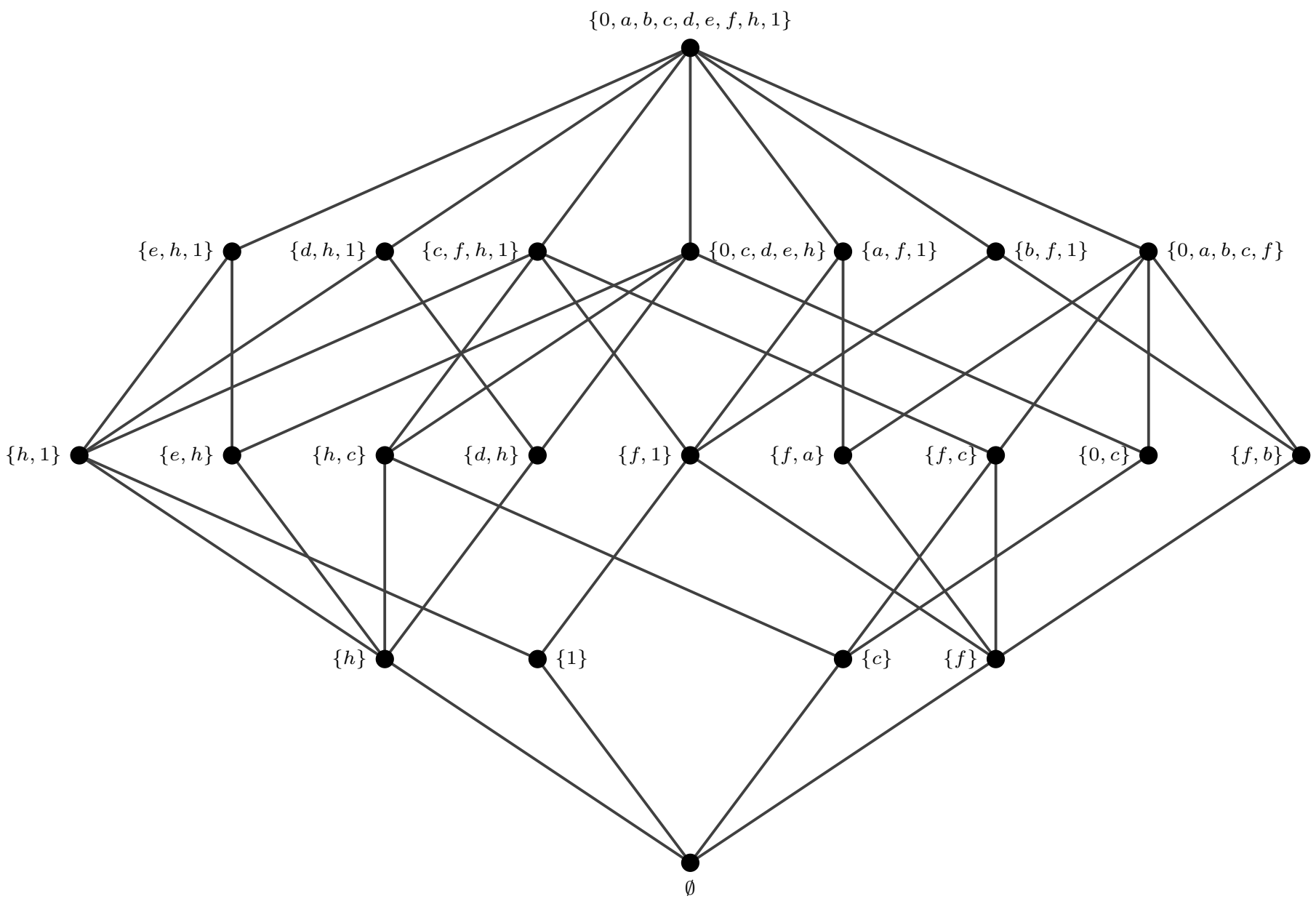
\begin{figure}[th]
		\begin{center}
			\resizebox{\textwidth}{!}{
				\begin{tikzpicture}[scale=1.3]
					\Vertex[size=0.01,x=0,y=0,label=${\emptyset}$,position=below,style={color=black}]{0}
					\Vertex[size=0.01,x=-3,y=2,label=${\{h\}}$,position=left,style={color=black}]{a}
					\Vertex[size=0.01,x=3,y=2,label=${\{f\}}$,position=left,style={color=black}]{b}
					\Vertex[size=0.1,x=-1.5,y=2,label=${\{1\}}$,position=right,style={color=black}]{c}
					\Vertex[size=0.01,x=1.5,y=2,label=${\{c\}}$,position=right,style={color=black}]{d}
					\Vertex[size=0.01,x=-6,y=4,label=${\{h,1\}}$,position=left,style={color=black}]{e}
					\Vertex[size=0.01,x=0,y=4,label=${\{f,1\}}$,position=left,style={color=black}]{f}
					\Vertex[size=0.01,x=-3,y=4,label=${\{h,c\}}$,position=left,style={color=black}]{g}
					\Vertex[size=0.01,x=-1.5,y=4,label=${\{d,h\}}$,position=left,style={color=black}]{h}
					\Vertex[size=0.01,x=-4.5,y=4,label=${\{e,h\}}$,position=left,style={color=black}]{i}
					\Vertex[size=0.01,x=1.5,y=4,label=${\{f,a\}}$,position=left,style={color=black}]{j}
					\Vertex[size=0.01,x=3,y=4,label=${\{f,c\}}$,position=left,style={color=black}]{k}
					\Vertex[size=0.01,x=6,y=4,label=${\{f,b\}}$,position=left,style={color=black}]{l}
					\Vertex[size=0.01,x=4.5,y=4,label=${\{0,c\}}$,position=left,style={color=black}]{m}
					\Vertex[size=0.01,x=-4.5,y=6,label=${\{e,h,1\}}$,position=left,style={color=black}]{n}
					\Vertex[size=0.01,x=-3,y=6,label=${\{d,h,1\}}$,position=left,style={color=black}]{p}
					\Vertex[size=0.01,x=-1.5,y=6,label=${\{c,f,h,1\}}$,position=left,style={color=black}]{q}
					\Vertex[size=0.01,x=1.5,y=6,label=${\{a,f,1\}}$,position=right,style={color=black}]{r}
					\Vertex[size=0.01,x=3,y=6,label=${\{b,f,1\}}$,position=right,style={color=black}]{s}
					\Vertex[size=0.01,x=0,y=6,label=${\{0,c,d,e,h\}}$,position=right,style={color=black}]{t}
					\Vertex[size=0.01,x=4.5,y=6,label=${\{0,a,b,c,f\}}$,position=right,style={color=black}]{u}
					\Vertex[size=0.01,x=0,y=8,label=${\{0,a,b,c,d,e,f,h,1\}}$,position=above,style={color=black}]{1}

					\Edge[lw=1pt](0)(a)
					\Edge[lw=1pt](0)(b)
					\Edge[lw=1pt](0)(c)
					\Edge[lw=1pt](0)(d)
					
					\Edge[lw=1pt](a)(e)
					\Edge[lw=1pt](a)(h)
					\Edge[lw=1pt](a)(g)
					\Edge[lw=1pt](a)(i)
					
					\Edge[lw=1pt](b)(f)
					\Edge[lw=1pt](b)(j)
					\Edge[lw=1pt](b)(k)
					\Edge[lw=1pt](b)(l)

					\Edge[lw=1pt](c)(f)
					\Edge[lw=1pt](c)(e)
					
					\Edge[lw=1pt](d)(g)
					\Edge[lw=1pt](d)(k)
					\Edge[lw=1pt](d)(m)

					\Edge[lw=1pt](e)(n)
					\Edge[lw=1pt](e)(p)
					\Edge[lw=1pt](e)(q)
					
					\Edge[lw=1pt](f)(q)
					\Edge[lw=1pt](f)(r)
					\Edge[lw=1pt](f)(s)
					
					\Edge[lw=1pt](g)(q)
					\Edge[lw=1pt](g)(t)
					
					\Edge[lw=1pt](h)(p)
					\Edge[lw=1pt](h)(t)
					
					\Edge[lw=1pt](i)(n)
					\Edge[lw=1pt](i)(t)
					
					\Edge[lw=1pt](j)(r)
					\Edge[lw=1pt](j)(u)
					
					\Edge[lw=1pt](k)(q)
					\Edge[lw=1pt](k)(u)
					
					\Edge[lw=1pt](l)(s)
					\Edge[lw=1pt](l)(u)
					
					\Edge[lw=1pt](m)(u)
					\Edge[lw=1pt](m)(t)
					
					\Edge[lw=1pt](1)(n)
					\Edge[lw=1pt](1)(p)
					\Edge[lw=1pt](1)(q)
					\Edge[lw=1pt](1)(r)
					\Edge[lw=1pt](1)(s)
					\Edge[lw=1pt](1)(t)
					\Edge[lw=1pt](1)(u)
					
				\end{tikzpicture}
			}
			\caption{Lattice $(\CS(\mathfrak{L}) \cup \{\emptyset\}, \subseteq)$ for the hyperlattice considered in Example~\ref{eg1} }\label{fig3}
		\end{center}
	\end{figure}
\end{example}
\section{The set $\CS(\mathfrak{L})$ with respect to orders $\leq_b$ and $\sqsubseteq$}
In this section, we investigate the set $\CS(\mathfrak{L})$ with respect to two new orders, $\leq_b$ and~$\sqsubseteq$. We determine classes of hyperlattices where the two orders are equivalent and give examples to show that they do not coincide in general. We explore the relationship between $\CS(\mathfrak{L})$ and $I(\mathfrak{L}) \times F(\mathfrak{L})$ in $P$-hyperlattices. 
\begin{proposition}
Let $\leq_b$ be a relation defined on $\CS(\mathfrak{L})$ as follows: for $C_1, C_2 \in \CS(\mathfrak{L})$,
	$C_1 \leq_b C_2$ if and only if for every $a \in C_1$, there exists $b \in C_2$ such that $a \preceq b$ and for every $b \in C_2$ there exists $a \in C_1$ such that $a \preceq b$. Then $(\CS(\mathfrak{L}), \leq_b)$ forms a poset. 
\end{proposition}
\begin{proof}
	Clearly $\leq_b$ is reflexive and transitive. Now suppose $C_1 \leq_b C_2$ and $C_2 \leq_b C_1$. Let $a \in C_1$. Then there exists $b \in C_2$ such that $a \preceq b$. Since $C_2 \leq_b C_1$, there exists $b' \in C_2$ such that $b' \preceq a$. Hence $b' \preceq a \preceq b$ where $b,b' \in C_2$. By the convexity of $C_2$, it follows that $a \in C_2$. Thus $C_1 \subseteq C_2$. A similar argument gives $C_2 \subseteq C_1$, and hence $C_1=C_2$. Therefore $\leq_b$ is anti-symmetric. Hence, $(\CS(\mathfrak{L}), \leq_b)$ is a poset.
	
\end{proof}
\begin{proposition}\label{proposition:idealorder}
	Let $\mathfrak{L}$ be a $P$-hyperlattice. The relation $\sqsubseteq$ defined on $\CS(\mathfrak{L})$ by: for $C_1, C_2 \in \CS(\mathfrak{L})$, $C_1 \sqsubseteq C_2$ if and only if $\langle C_1 ] \subseteq \langle C_2]$ and $[C_1\rangle \supseteq [C_2\rangle$, is a partial order on $\CS(\mathfrak{L})$.
\end{proposition} 

\begin{proof}
	By definition, $\sqsubseteq$ is reflexive and transitive. 
	
	Now suppose $C_1 \sqsubseteq C_2$ and $C_2 \sqsubseteq C_1$. Then $\langle C_1 ] \subseteq \langle C_2]$, $[C_1\rangle \supseteq [C_2\rangle$ and $\langle C_2 ] \subseteq \langle C_1]$, $[C_2\rangle \supseteq [C_1 \rangle$. Thus $\langle C_1 ] = \langle C_2]$ and $[C_1\rangle = [C_2\rangle$. Since $\mathfrak{L}$ is a $P$-hyperlattice, using Theorem~\ref{thmP}, we get $C_1=\langle C_1 ] \cap [C_1\rangle$ and $C_2 = \langle C_2] \cap [C_2 \rangle$ from which it follows that $C_1=C_2$. Therefore $\sqsubseteq$ is anti-symmetric. 
	
\end{proof}
\begin{proposition}\label{proposition:bidomimpliesidealorder}
	Let $A, B  \in \CS(\mathfrak{L})$. If $A \leq_b B$ then $A \sqsubseteq B$.
\end{proposition}
\begin{proof}
	Let $A, B \in \CS(\mathfrak{L})$. Suppose $A \leq_b B$. Let $x \in A$. Then, there exists $z \in B$ such that $x \preceq z$. Since $z \in B \subseteq\langle B]$, and $x \preceq z$, we get $x \in \langle B]$. Therefore, $A \subseteq \langle B]$ which implies, $ \langle A] \subseteq \langle B]$. Now, let $x \in B$. Then, there exists $y \in A$ such that $y \preceq x$. Since $y \in A \subseteq [A\rangle $ and $y \preceq x$, we get that $x \in [A\rangle $. Hence $B \subseteq [A \rangle$ which implies $[B \rangle \subseteq [A\rangle $. Thus, $A \sqsubseteq B$.
\end{proof}

As illustrated in Example~\ref{eg10}, the relation $\sqsubseteq$ is not a partial order on an arbitrary strong join hyperlattice. In addition, Example~\ref{eg10} provides a counterexample to the equivalence of $\leq_b$ and $\sqsubseteq$.
\begin{example}\label{eg10}
	Consider the hyperlattice $(\mathfrak{L},\hjoin,\wedge)$ where $\hjoin$ and $\wedge$ are given in Table~\ref{Table6}. The hyperideals, type-III hyperfilters and convex subhyperlattices are listed in Table~\ref{Table8}. The poset $(\CS(\mathfrak{L}),\leq_b)$ is given in Figure~\ref{fig14}. The relation $\sqsubseteq$ is not a partial order on $\CS(\mathfrak{L})$ because the sets $\{x_1\}$ and $\{x_1,x_2,1\}$ violate antisymmetry. Further, $\{x_1,x_2,1\} \sqsubseteq \{x_1\}$, we have $\{x_1,x_2,1\} \nleq_B \{x_1\}$  demonstrating that $\sqsubseteq$ does not imply $\leq_b$. 
	
	\begin{table}[ht]

				\begin{tabular}{|c|c|c|c|c|}
					\hline 
					$\hjoin$ &$ 0 $ & $ x_1 $ & $ x_2 $ & $ 1 $  \\ 
					\hline 
					$0$ &$ \{0\}  $ & $ \{x_1,x_2,1\} $ & $\{x_2,1\}  $ & $\{1\}$  \\  \hline
					$x_1$ &$ \{x_1,x_2,1\} $ & $  \{x_1\}$ & $\{x_2,1\}  $ & $\{1\}$  \\  \hline
					$x_2$ &$ \{x_2,1\} $ & $\{x_2,1\}  $ & $ \{x_2\} $ & $\{1\}$ \\  \hline
					$1$ &$ \{1\} $ & $ \{1\} $ & $\{1\}  $ & $ \{1\} $ \\  \hline
				\end{tabular}
				\qquad
				\begin{tabular}{|c|c|c|c|c|}
					\hline 
					$\wedge$ &$ 0 $ & $x_1  $ & $x_2  $ & $1  $ \\  \hline
					$0$ &$ 0 $ & $ 0 $ & $0  $ & $0  $  \\  \hline
					$x_1$ &$0  $ & $x_1  $ & $x_1  $ & $x_1  $ \\  \hline
					$x_2$ &$0  $ & $x_1  $ & $x_2  $ & $x_2  $  \\  \hline
					$1$ &$0  $ & $x_1  $ & $x_2  $ & $1  $  \\  \hline		
				\end{tabular} 
			\caption{$\hjoin$ and $\wedge$ on $\mathfrak{L}$ considered in Example~\ref{eg10}} \label{Table6}

	\end{table}

		\begin{table}[ht]
				\begin{tabular}{|c|c|cc|}
					\hline
					{Hyperideals}&{ Type-III Hyperfilters}&\multicolumn{2}{c|}{{Convex subhyperlattices}}\\
					\hline
					$I_1=\{0\}$ &$F_1=\{1\}$ &$C_1=\{0\}$&$C_{6}=\{x_1,x_2,1\}$\\
					$I_2=\{0,x_1,x_2,1\}$ &	$F_2=\{x_2,1\}$&$C_2=\{x_1\}$ & $C_{7}=\{0,x_1,x_2,1\}$\\
					& $F_3=\{x_1,x_2,1\}$ &$C_3=\{x_2\}$& \\
					$I_4=L$& $F_=\{0,x_1,x_2,1\}$ &$C_4=\{1\}$	& \\
					&&$C_5=\{x_2,1\}$ &\\
					
					\hline
				\end{tabular}
			\caption{Hyperideals, Type-III Hyperfilters and Convex subhyperlattices of $(\mathfrak{L},\hjoin,\wedge)$ considered in Example~\ref{eg10}}\label{Table8}

	\end{table}
	
	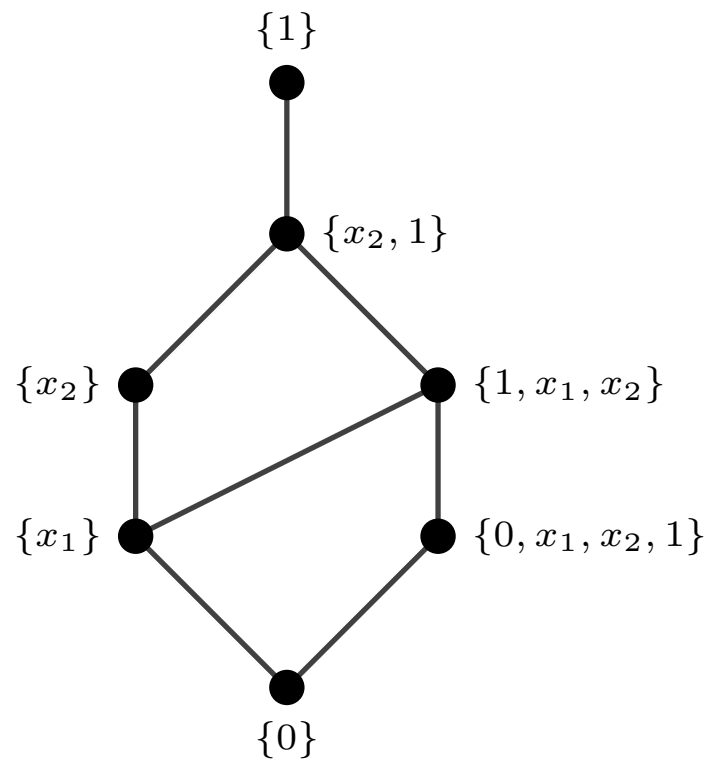
\begin{figure}[th]
		\begin{center}
			
			\begin{tikzpicture}
				\Vertex[size=0.01,x=0,y=0.5,label=${\{0\}}$,position=below,style={color=black}]{0}
				
				\Vertex[size=0.01,x=-1,y=1.5,label=${\{x_1\}}$,position=left,style={color=black}]{b}
				\Vertex[size=0.1,x=1,y=1.5,label=${\{0,x_1,x_2,1\}}$,position=right,style={color=black}]{a}
				\Vertex[size=0.01,x=1,y=2.5,label=${\{1,x_1,x_2\}}$,position=right,style={color=black}]{c}
				
				\Vertex[size=0.01,x=0,y=3.5,label=${\{x_2,1\}}$,position=right,style={color=black}]{d}
				
				\Vertex[size=0.01,x=-1,y=2.5,label=${\{x_2\}}$,position=left,style={color=black}]{e}
				\Vertex[size=0.01,x=0,y=4.5,label=${\{1\}}$,position=above,style={color=black}]{f}

				\Edge[lw=1pt](0)(a)
				\Edge[lw=1pt](0)(b)
				
				\Edge[lw=1pt](a)(c)
				\Edge[lw=1pt](c)(d)
				
				\Edge[lw=1pt](b)(e)
				\Edge[lw=1pt](e)(d)
				
				\Edge[lw=1pt](b)(c)
				\Edge[lw=1pt](f)(d)
				
			\end{tikzpicture}
			
			\caption{Poset $(\CS(\mathfrak{L}) , \leq_b)$ for the hyperlattice considered in Example~\ref{eg10}  }\label{fig14}
		\end{center}
	\end{figure}
\end{example}
\begin{proposition}[\cite{Nayak2025}] \label{idealofPHL}
	Let $ H \subseteq \mathfrak{L}$ be non-empty, where $\mathfrak{L}$ is a $P$-hyperlattice. Then 
	\begin{enumerate}[label=\roman*]
		\item $\langle H]=\left\{x \in \mathfrak{L} : x \preceq y ~\text{for some}~ y \in h_1 \hjoin\limits^P h_2 \hjoin\limits^P \dots\hjoin\limits^P h_n, \right.$
		\\
		$\left.\phantom{\langle H]=\{x\in \mathfrak L\hjoin\limits^P} h_i \in H, n \in \mathbb{N}, 1 \leq i \leq n\right \}.$
		\item $[H \rangle =\left \{x \in \mathfrak{L}:  h_1 \wedge h_2 \wedge \dots \wedge h_n \preceq x , h_i \in H, n \in \mathbb{N}, 1 \le i \le n \right \}.$
	\end{enumerate} 
\end{proposition}

\noindent The following proposition shows that orders $\leq_b$ and $\sqsubseteq$ are equivalent in a $P$-hyperlattice.
\begin{proposition}\label{twoorderequivalent}
	Let $\mathfrak{L}$ be a $P$-hyperlattice. Then the orders $\leq_b$ and $\sqsubseteq$ are equivalent.
\end{proposition}

\begin{proof}
	Assume that $A \sqsubseteq B$ where $A,B \in \CS(\mathfrak{L})$. Let $x \in A$. Then since $\langle A ] \subseteq \langle B]$, we get $x \in \langle B]$. This implies $x \preceq y$ for some $y \in h_1 \hjoin\limits^P h_2 \hjoin\limits^P \dots \hjoin\limits^P h_n$ where $h_i \in B$ ($1 \leq i \leq n$) and $n \in \mathbb{N}$. Since $B$ is a subhyperlattice, we get $ y \in B$. For the verification of the second condition, let $z \in B$. Then since $B \subseteq [B\rangle \subseteq [A\rangle $, we have $g_1 \wedge g_2 \wedge \dots \wedge g_k \preceq z$ where $g_i \in A$ ($1 \leq i \leq k$) and $k \in \mathbb{N}$. Now the fact that $A$ is a subhyperlattice gives that $w=g_1 \wedge g_2 \wedge \dots \wedge g_k \in A$. Thus, we have a $w \in A$ such that $w \preceq z$. Therefore, $A \leq_b B$. Hence, the orders $\leq_b$ and $\sqsubseteq$ are equivalent by Proposition~\ref{proposition:bidomimpliesidealorder}.
\end{proof}

\begin{example}\label{eg7}
	Let $\mathfrak{L}=\{a,b,c,d,e\}$. The hyperoperation $\hjoin$ and classical operation $\wedge$ are given in Table~\ref{Table5}. Then $(\mathfrak{L},\hjoin,\wedge)$ forms a join hyperlattice but it is not a strong hyperlattice. In this hyperlattice, the orders $\leq_b$ and $\sqsubseteq$ on $\CS(\mathfrak{L})$ are not equivalent. Here, $(\CS(\mathfrak{L}), \leq_b)$ is a poset (infact a lattice) which is shown in Figure~\ref{fig11}. But $\sqsubseteq$ is not a partial order on $\CS(\mathfrak{L})$. If we take $A=\{b\}$ and $B=\{a,b\}$, then we get $A \sqsubseteq B$ and $B \sqsubseteq A$. But $A \neq B$. 
	
	\begin{table}[ht]
	
				\begin{tabular}{|c|c|c|c|c|c|}
					\hline 
					$\hjoin$ &$ a $ & $ b $ & $ c $ & $ d $&$e$  \\ 
					\hline 
					$a$ &$ \{a\}  $ & $ \{a,b\} $ & $\{a,c\}  $ & $\{a,d\}  $& $\{a,e\}$  \\  \hline
					$b$ &$ \{a,b\} $ & $  \{b\}$ & $\{b,c,e\}  $ & $\{b,d,e\}  $& $\{b,e\}$  \\  \hline
					$c$ &$ \{a,c\} $ & $\{b,c,e\}  $ & $ \{c\} $ & $\{c,d\}  $& $\{c,e\}$  \\  \hline
					$d$ &$ \{a,d\} $ & $ \{b,d,e\} $ & $\{c,d\}  $ & $ \{d\} $& $\{d,e\}$  \\  \hline
					$e$ &$ \{a,e\} $ & $ \{b,e\} $ & $\{c,e\}  $ & $\{d,e\}  $& $\{e\}$  \\  \hline
				\end{tabular}
				\qquad
				\begin{tabular}{|c|c|c|c|c|c|}
					\hline 
					$\wedge$ &$ a $ & $b  $ & $c  $ & $d  $& $e$  \\  \hline
					$a$ &$ a $ & $ a $ & $a  $ & $a  $& $a$  \\  \hline
					$b$ &$a  $ & $b  $ & $a  $ & $a  $& $b$  \\  \hline
					$c$ &$a  $ & $a  $ & $c  $ & $c  $& $c$  \\  \hline
					$d$ &$a  $ & $a  $ & $c  $ & $d  $& $d$  \\  \hline
					$e$ &$ a $ & $b  $ & $c  $ & $d  $& $e$  \\  \hline
					
				\end{tabular}
		\caption{$\hjoin$ and $\wedge$ on $\mathfrak{L}$ considered in Example~\ref{eg7} } \label{Table5}
	\end{table}
	
	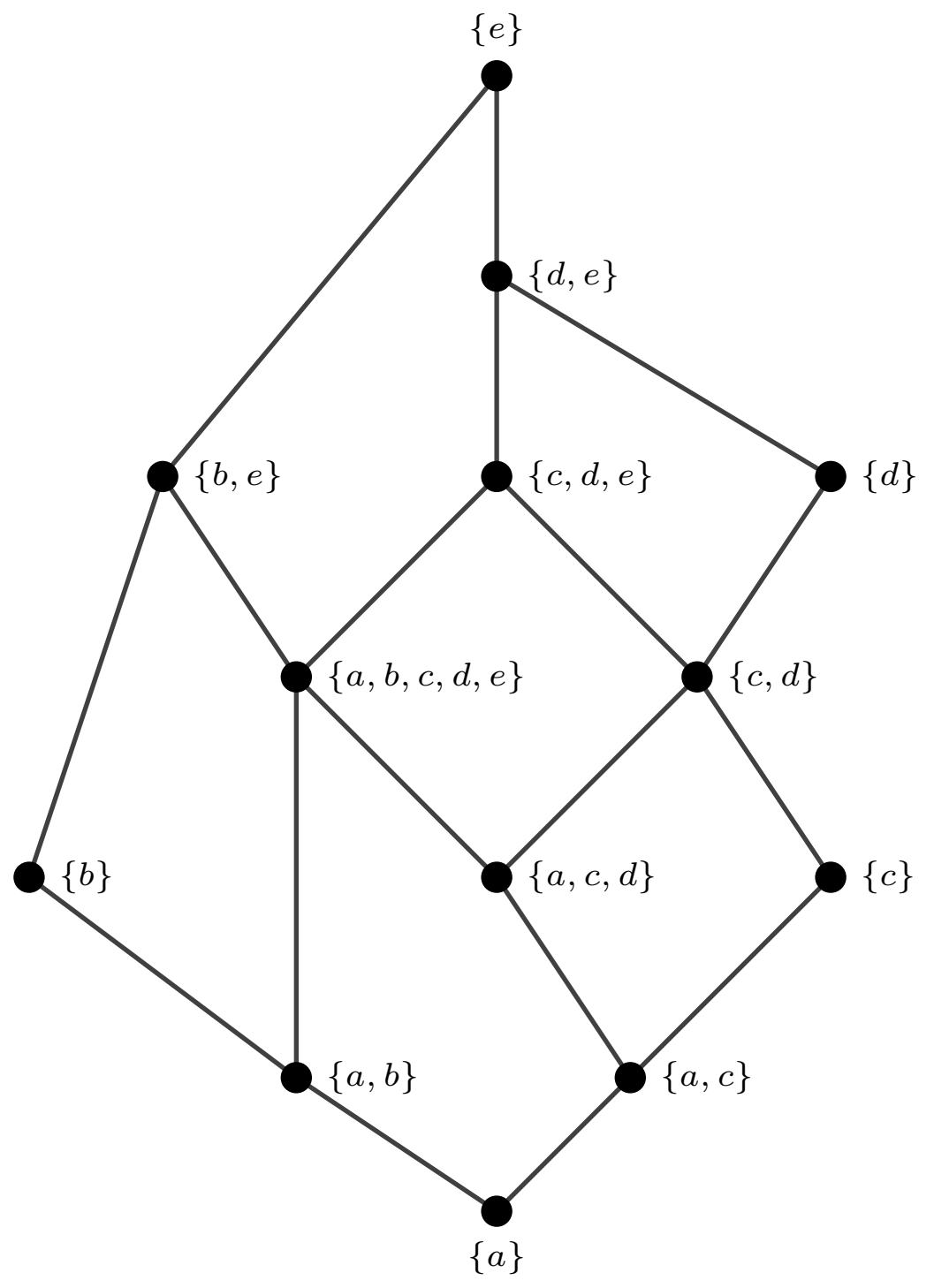
\begin{figure}[th]
		\begin{center}
			
			\begin{tikzpicture}
				\Vertex[size=0.01,x=0,y=0.5,label=${\{a\}}$,position=below,style={color=black}]{0}
				
				\Vertex[size=0.01,x=-1.5,y=1.5,label=${\{a,b\}}$,position=right,style={color=black}]{a}
				\Vertex[size=0.1,x=1,y=1.5,label=${\{a,c\}}$,position=right,style={color=black}]{b}
				\Vertex[size=0.01,x=-3.5,y=3,label=${\{b\}}$,position=right,style={color=black}]{d}
				
				\Vertex[size=0.01,x=0,y=3,label=${\{a,c,d\}}$,position=right,style={color=black}]{c}
				
				\Vertex[size=0.01,x=2.5,y=3,label=${\{c\}}$,position=right,style={color=black}]{e}
				
				\Vertex[size=0.01,x=-1.5,y=4.5,label=${\{a,b,c,d,e\}}$,position=right,style={color=black}]{f}
				
				\Vertex[size=0.01,x=1.5,y=4.5,label=${\{c,d\}}$,position=right,style={color=black}]{g}
				\Vertex[size=0.01,x=-2.5,y=6,label=${\{b,e\}}$,position=right,style={color=black}]{h}
				\Vertex[size=0.01,x=0,y=6,label=${\{c,d,e\}}$,position=right,style={color=black}]{i}
				\Vertex[size=0.01,x=2.5,y=6,label=${\{d\}}$,position=right,style={color=black}]{j}
				\Vertex[size=0.01,x=0,y=7.5,label=${\{d,e\}}$,position=right,style={color=black}]{k}
				\Vertex[size=0.01,x=0,y=9,label=${\{e\}}$,position=above,style={color=black}]{1}
				
				\Edge[lw=1pt](0)(a)
				\Edge[lw=1pt](0)(b)
				
				\Edge[lw=1pt](a)(d)
				\Edge[lw=1pt](a)(f)
				
				\Edge[lw=1pt](b)(c)
				\Edge[lw=1pt](b)(e)
				
				\Edge[lw=1pt](c)(f)
				\Edge[lw=1pt](c)(g)
				
				\Edge[lw=1pt](d)(h)
				
				\Edge[lw=1pt](e)(g)
				
				\Edge[lw=1pt](f)(h)
				\Edge[lw=1pt](f)(i)
				
				\Edge[lw=1pt](g)(i)
				\Edge[lw=1pt](g)(j)
				
				\Edge[lw=1pt](h)(1)
				
				\Edge[lw=1pt](i)(k)
				
				\Edge[lw=1pt](j)(k)
				
				\Edge[lw=1pt](k)(1)

			\end{tikzpicture}
			
			\caption{Lattice $(\CS(\mathfrak{L}) , \leq_b)$ for the hyperlattice considered in Example~\ref{eg7}  }\label{fig11}
		\end{center}
	\end{figure}
\end{example}

\begin{remark}\label{rem:idealbyCSH}
	Let $\mathfrak{L}$ be a $P$-hyperlattice. Then,
	\begin{enumerate}[label=\roman*]
		\item for a non-empty subset $A$ of $\mathfrak{L}$, we have $\Big \langle \langle A\rangle \Big ]= \langle A]$ and $\Big[ \langle A \rangle \Big \rangle =[A\rangle$.
		\item for two non-empty subsets $A$ and $B$ of $\mathfrak{L}$, suppose that for every $a \in A$, there is a $b \in B$ such that $a \preceq b$ and for every $b \in B$, there is a $a \in A$ such that $a \preceq b$. Then $A \subseteq \langle B]$ and $B \subseteq [A\rangle$. Then by Proposition~\ref{twoorderequivalent}, we get $\langle A \rangle \leq_b \langle B \rangle$.
		\item for $a,b \in \mathfrak{L}$ such that $a \preceq b$, we have $\langle a \rangle \leq_b \langle b \rangle$.
		\item The map $f: \mathfrak{L} \rightarrow (\CS(\mathfrak{L}),\leq_b)$ defined by $f(a)=\langle a \rangle$ is order preserving. 
		
	\end{enumerate}
\end{remark}

\begin{proposition}\label{idealorderlattice}
	Let $\mathfrak{L}$ be a $P$-hyperlattice. Then the poset $(\CS(\mathfrak{L}),\sqsubseteq)$ forms a lattice with the following meet and join operations: For $C_1,C_2 \in \CS(\mathfrak{L})$,
	\begin{eqnarray*}
		C_1 \wedge C_2 &=& \left(\langle C_1] \wedge \langle C_2] ) \cap ([C_1\rangle \wedge [C_2 \rangle\right)\\
		C_1 \vee C_2 &=&(\langle C_1] \vee \langle C_2]) \cap ([C_1\rangle \vee [C_2\rangle )
	\end{eqnarray*}
	where $\langle C_1] \wedge \langle C_2]=\langle C_1] \cap \langle C_2]$,  $\langle C_1] \vee \langle C_2] =\Big \langle \langle C_1] \cup \langle C_2]\Big]$ are the meet and join operations in $I(\mathfrak{L})$, and $[C_1\rangle \wedge [C_2 \rangle =\Big[[C_1\rangle \cup [C_2 \rangle \Big \rangle$, $[C_1\rangle \vee [C_2 \rangle = [C_1\rangle \cap [C_2 \rangle$ are the meet and join operations in $F(\mathfrak{L})$.
\end{proposition}
\begin{proof}
	Let $M= \left(\langle C_1] \wedge \langle C_2] ) \cap ([C_1\rangle \wedge [C_2 \rangle\right)$ and $N=\left(\langle C_1] \vee \langle C_2]) \cap ([C_1\rangle \vee [C_2\rangle \right)$. First we show that $M$ and $N$ are non-empty. Choose $x \in C_1$ and $y \in C_2$. Then $x \wedge y \in \langle C_1], \langle C_2]$. Also, $x \in C_1 \subseteq [C_1\rangle  \subseteq \Big[ [C_1\rangle \cup [ C_2 \rangle \Big\rangle$. Similarly, $y \in \Big[ [C_1\rangle \cup [ C_2 \rangle \Big\rangle$. By the definition of hyperfilter, we get that $x \wedge y \in \Big[ [C_1\rangle \cup [ C_2 \rangle \Big\rangle$. Thus $M$ is non-empty. To show $N$ is non-empty, let $a \in C_1 \subseteq [C_1\rangle $ and $b \in C_2 \subseteq [C_2\rangle$. This implies $a \hjoin\limits^P b\subseteq  [C_1\rangle \cap [C_2\rangle$ by the definition of hyperfilter. We also have $a, b \in \langle C_1 ] \cup \langle C_2 ] \subseteq \Big\langle\langle C_1] \cup \langle C_2] \Big]$. Since $\Big\langle\langle C_1] \cup \langle C_2] \Big]$ is a hyperideal, we get $a \hjoin\limits^P b  \subseteq\Big\langle\langle C_1] \cup \langle C_2] \Big]$. This implies that $a \hjoin b \subseteq N$, and hence, $N$ is non-empty. Thus by Theorem~\ref{thmC}, it follows that $M$ and $N$ are convex subhyperlattices. Since $\mathfrak{L}$ is a $P$-hyperlattice, $M$ attains a unique representation in terms of $\langle M]$ and $[M\rangle$. By the uniqueness, we get that $\langle M] =\langle C_1 ] \cap \langle C_2]$ and $[M\rangle = \Big[ [C_1\rangle \cup [ C_2 \rangle \Big\rangle$. By a similar argument, we can show that $\langle N] =\Big \langle \langle C_1 ] \cup \langle C_2] \Big]$ and $[N \rangle = [C_1\rangle \cap [ C_2 \rangle$. Using these facts, we conclude that $M$ and $N$ are greatest lower bound and least upper bound of $\{C_1, C_2\}$ respectively, in $(\CS(\mathfrak{L}),\sqsubseteq)$. 
	
\end{proof}

\begin{corollary}
	$(\CS(\mathfrak{L}), \leq_b)$ is a lattice for a $P$-hyperlattice $\mathfrak{L}$.
\end{corollary}
\begin{proof}
	Follows by Proposition~\ref{twoorderequivalent} and Proposition~\ref{idealorderlattice}.
\end{proof}
The following example shows that $(\CS(\mathfrak{L}),\leq_b)$ need not be a lattice for an arbitrary join hyperlattice.
\begin{example}\label{eg8}
	For $\mathfrak{L}=\{a_1,b_1,c_1\}$, let the operation $\wedge$ and hyperoperation $\hjoin$ be as given in Table~\ref{Table7}. Then the poset $(\CS(\mathfrak{L}),\leq_b)$ is given in Figure~\ref{fig12}. Here $(\CS(\mathfrak{L}),\leq_b)$ is not a lattice.
	
	\begin{table}[ht] 
				\begin{tabular}{|c|c|c|c|}
					\hline 
					$\hjoin$ &$ a_1 $ & $ b_1 $ & $ c_1 $  \\ 
					\hline 
					$a_1$ &$ \{a_1\}  $ & $ \{b_1\} $ & $\{c_1\}  $  \\  \hline
					$b_1$ &$ \{b_1\} $ & $  \{b_1\}$ & $\{b_1,c_1\}  $   \\  \hline
					$c_1$ &$ \{c_1\} $ & $\{b_1,c_1\}  $ & $ \{c_1\} $  \\  \hline
				\end{tabular}
				\qquad
				\begin{tabular}{|c|c|c|c|}
					\hline 
					$\wedge$ &$ a_1 $ & $b_1  $ & $c_1  $  \\  \hline
					$a_1$ &$ a_1 $ & $ a_1 $ & $a_1  $  \\  \hline
					$b_1$ &$a_1  $ & $b_1  $ & $a_1  $   \\  \hline
					$c_1$ &$a_1  $ & $a_1  $ & $c_1  $  \\  \hline
				\end{tabular}
		\caption{$\hjoin$ and $\wedge$ on $\mathfrak{L}$ considered in Example~\ref{eg8}}  \label{Table7}

	\end{table}
	
	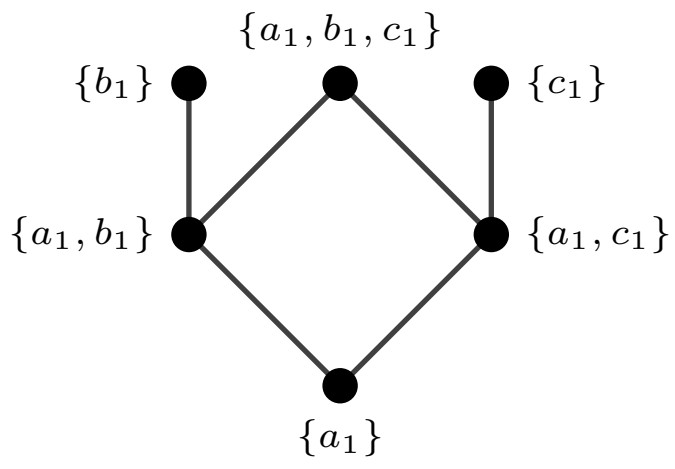
\begin{figure}[th]
		\begin{center}
			
			\begin{tikzpicture}
				\Vertex[size=0.01,x=0,y=0.5,label=${\{a_1\}}$,position=below,style={color=black}]{0}
				
				\Vertex[size=0.01,x=-1,y=1.5,label=${\{a_1,b_1\}}$,position=left,style={color=black}]{b}
				\Vertex[size=0.1,x=1,y=1.5,label=${\{a_1,c_1\}}$,position=right,style={color=black}]{a}
				\Vertex[size=0.01,x=1,y=2.5,label=${\{c_1\}}$,position=right,style={color=black}]{c}
				
				\Vertex[size=0.01,x=0,y=2.5,label=${\{a_1,b_1,c_1\}}$,position=above,style={color=black}]{d}
				
				\Vertex[size=0.01,x=-1,y=2.5,label=${\{b_1\}}$,position=left,style={color=black}]{e}

				\Edge[lw=1pt](0)(a)
				\Edge[lw=1pt](0)(b)
				
				\Edge[lw=1pt](a)(c)
				\Edge[lw=1pt](a)(d)
				
				\Edge[lw=1pt](b)(e)
				\Edge[lw=1pt](b)(d)

			\end{tikzpicture}
			
			\caption{Poset $(\CS(\mathfrak{L}) , \leq_b)$ for the hyperlattice considered in Example~\ref{eg8}  }\label{fig12}
		\end{center}
	\end{figure}
\end{example}
\begin{definition}
	Let $X$ be a non-empty subset of $\mathfrak{L}$.
	\begin{enumerate}[label=\roman*]
		\item The initial segment of $X$ with respect to $\preceq$ is defined by
		\begin{equation*}
			\big \downarrow_\preceq(X)=\{y \in \mathfrak{L} \mid y \preceq x~\text{for some}~ x \in X\}
		\end{equation*}
		If $X=\{x\}$, we use the symbol $\downarrow_\preceq (x)$ in place of $\downarrow_\preceq (\{x\})$.
		\item The final segment of $X$ with respect to $\preceq$ is defined by
		\begin{equation*}
			\big \uparrow_\preceq(X)=\{y \in \mathfrak{L} \mid x \preceq y~\text{for some}~ x \in X\}
		\end{equation*}
		If $X=\{x\}$, we use the symbol $\uparrow_\preceq(x)$ in place of $\uparrow_\preceq (\{x\})$.
	\end{enumerate}
\end{definition}

\noindent Note that $X \subseteq  {\downarrow_\preceq(X)} \subseteq \langle X]$ and $X \subseteq {\uparrow_\preceq(X)} \subseteq [X \rangle $ for any non-empty subset $X$ of $\mathfrak{L}$.

\begin{example}\label{ch3:eg5}
	Consider the hyperlattice shown in Figure~\ref{fig8}. Suppose $X=\{e,f\}$. Then we get~$\downarrow_\preceq(X)=\{0,a,b,c,e,f\}$, $\uparrow_\preceq(X)=\{e,f,1\}, [X\rangle =\{c,e,f,1\}$ and $\langle X] =\mathfrak{L}$. This shows that $\downarrow_\preceq(X) \subset \langle X]$ and $\uparrow_\preceq(X) \subset [X\rangle$. The initial segment of $X$ is shown by dotted lines and the final segment of $X$ is shown by dashed lines in Figure~\ref{fig8}.
\end{example}
Proposition~\ref{prop:csh} gives a characterization of those non-empty subsets of a $P$-hyperlattice that are convex subhyperlattices.
\begin{proposition}\label{prop:csh}
	A non-empty set $C$ of a $P$-hyperlattice $\mathfrak{L}$ is a convex subhyperlattice if and only if $\downarrow_\preceq(C) = \langle C], \uparrow_\preceq(C)= [C \rangle $ and $C$ is a convex set.
\end{proposition}
\begin{proof}
	Let $\emptyset \neq C \subseteq \mathfrak{L}$ be a convex subhyperlattice. Clearly, $C$ is a convex set. Already, $\downarrow_\preceq(C)\subseteq \langle C]$. Let $a \in \langle C]$. Then $a \preceq y$ for some $y \in h_1 \hjoin\limits^P h_2 \hjoin\limits^P \dots\hjoin\limits^P h_n$ where $h_i \in C, n \in \mathbb{N}$ and $1 \leq i \leq n$. Note that $h_1 \hjoin\limits^P h_2 \hjoin\limits^P \dots\hjoin\limits^P h_n \subseteq C$ since $C$ is a subhyperlattice. Thus we have $a \preceq y$ where $y \in C$. This implies $a \in \downarrow_\preceq~(C)$ and hence $\langle C] =  \downarrow_\preceq(C)$. Now, let $b \in [C \rangle$. Then $g_1 \wedge g_2 \wedge \dots \wedge g_s \preceq b$ for some $g_i \in C, s \in \mathbb{N}$. Since $C$ is a subhyperlattice, $g_1 \wedge g_2 \wedge \dots \wedge g_s \in C$, so we get $b \in \uparrow_\preceq(C)$. This shows that $[C\rangle \subseteq \uparrow_\preceq(C)$ and hence, $[C\rangle = \uparrow_\preceq(C)$. Conversely, assume that $\downarrow_\preceq(C) = \langle C], \uparrow_\preceq(C)= [C \rangle $ and $C$ is a convex set. To show $C$ is a convex subhyperlattice, it remains to show that $C$ is a subhyperlattice. Clearly $C \subseteq \langle C] \cap [C\rangle$. Let $b \in \langle C] \cap [C\rangle=\downarrow_\preceq(C)  \cap \uparrow_\preceq(C)$. Then $x \preceq b \preceq y$ for some $x,y \in C$. Since $C$ is convex, we get $b \in C$. Thus $C= \langle C] \cap [C\rangle$ which implies $C$ is a subhyperlattice.
\end{proof}

\begin{remark}
	The convexity of $C$ is necessary in the converse part of Proposition~\ref{prop:csh}. In the $P$-hyperlattice shown in Figure~\ref{fig8}, consider the subhyperlattice $C=\{0,c,f\}$. Then $\uparrow_\preceq(C)=\mathfrak{L}=[C\rangle $ and $\downarrow_\preceq(C)=\{0,b,c,f\}=\langle C]$. Note that $C$ is not a convex set. 
	
	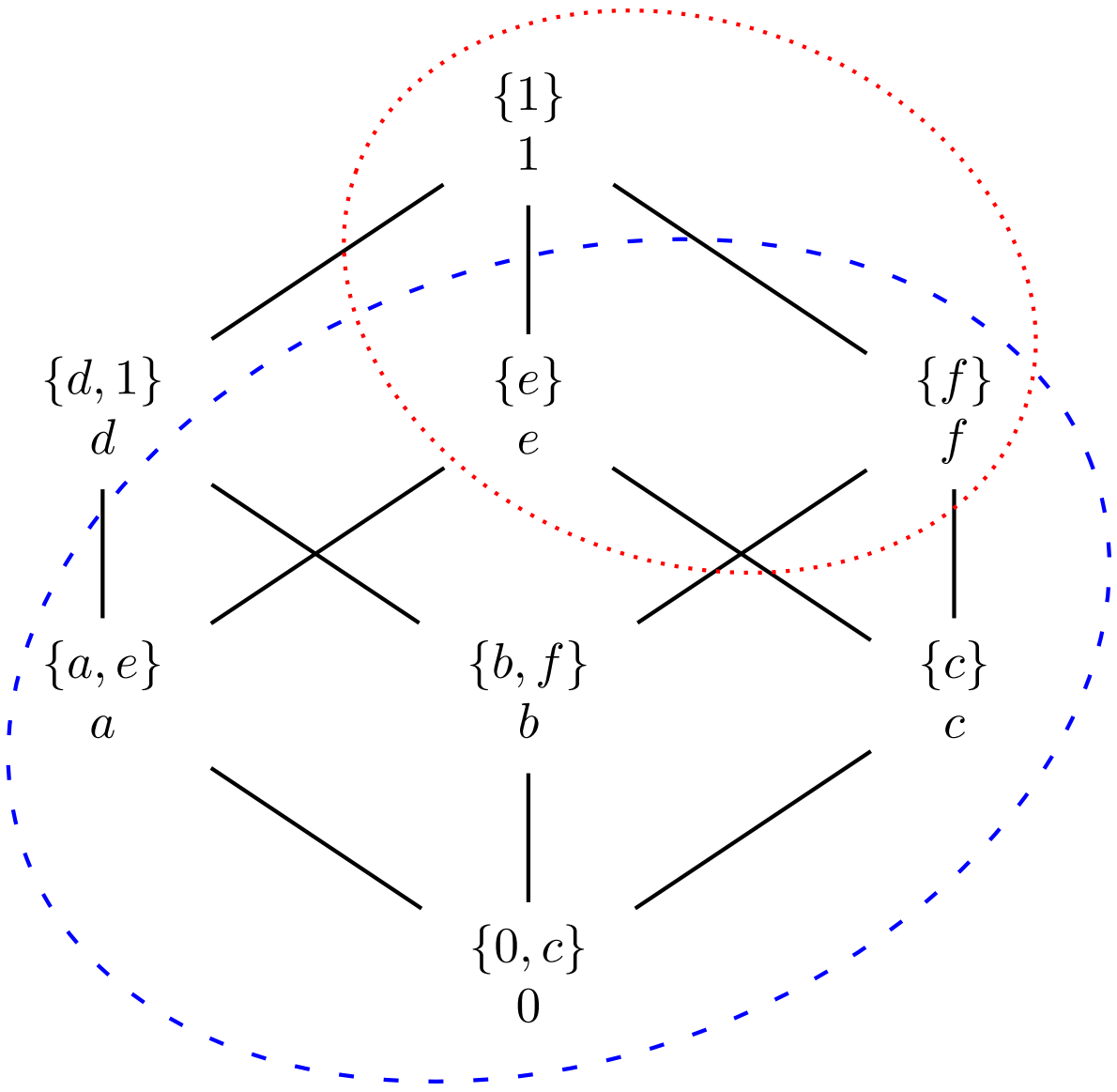
\begin{figure}[th]
		\begin{center}
			
			\begin{tikzpicture}
				
				\node  (h) at (0,4) {\begin{tabular}{c}
						$ \{1\} $\\$1 $
				\end{tabular}};
				
				\node  (d) at (-3,2) {\begin{tabular}{c}
						$ \{d,1\} $\\$d $
				\end{tabular}};;
				\node  (e) at (0,2) {\begin{tabular}{c}
						$ \{e\} $\\$e$
				\end{tabular}};
				\node  (f) at (3,2) {\begin{tabular}{c}
						$ \{f\} $\\$f $
				\end{tabular}};
				\node  (a) at (-3,0) {\begin{tabular}{c}
						$ \{a,e\} $\\$a $
				\end{tabular}};
				\node  (b) at (0,0) {\begin{tabular}{c}
						$ \{b,f\} $\\$b$
				\end{tabular}};
				\node  (c) at (3,0) {\begin{tabular}{c}
						$ \{c\} $\\$c$
				\end{tabular}};
				\node  (0) at (0,-2) {\begin{tabular}{c}
						$ \{0,c\} $\\$0 $
				\end{tabular}};
				\draw[thick]  (h)--(d)--(a)--(0)--(c)--(f)--(h);
				\draw[thick](d)--(b)--(f);
				\draw [thick](b)--(0);
				\draw[thick](a)--(e)--(c);
				\draw[thick](e)--(h);

				\draw [red,dotted,thick][rotate around={160:(1.3,3)}] (1.4,3.2) circle(2.5cm and 1.9cm);
				\draw [blue,loosely dashed,thick][rotate around={20:(-1,3)}] (-0.8,0) circle(4cm and 2.8cm);
				
			\end{tikzpicture}
			
			\caption{Hyperlattice considered in Example~\ref{ch3:eg5}}\label{fig8}
		\end{center}
	\end{figure}
\end{remark} 

\begin{proposition}[\cite{SoltaniLashkenari2015}]\label{prop: orderproperty}
	Let $(\mathfrak{L}, \hjoin, \wedge)$ be a strong join hyperlattice. For every $a,b \in \mathfrak{L}$, there exist $a_1,a_2 \in a \hjoin b$ such that $a \preceq a_1$ and $b \preceq a_2$.
\end{proposition}

\begin{proposition}\label{prop:4} 
	Let $X$ and $Y$ be two non-empty subsets of $(\mathfrak{L}, \hjoin\limits, \wedge)$, and \newline  $A=X \barwedge Y=\{x \wedge y\mid x \in X, y \in Y\}, ~B=X \hjoin Y=\mathlarger{\bigcup\limits}_{a \in X,b \in Y} a \hjoin b$, then
	\begin{enumerate}[label=\roman*]
		\item $\langle A] \subseteq \langle X] \cap \langle Y]$ and $[B\rangle \subseteq [X\rangle \cap [Y\rangle$,
		\item $[A\rangle = \Big[[X \rangle \cup [Y\rangle \Big \rangle$ and $\langle B] = \Big \langle \langle X] \cup \langle Y] \Big ]$,
	\end{enumerate}
\end{proposition}
\begin{proof}
	\begin{enumerate}[label=\roman*]
		\item Let $p \in A$. Then $p=x \wedge y$ for some $x \in X, y \in Y$. Since $x \wedge y \preceq x, y$, we get $p \in \langle X] \cap \langle Y]$. This gives $A \subseteq \langle X] \cap \langle Y]$, which implies $\langle A] \subseteq \langle X] \cap \langle Y]$. Now suppose $q \in B$. Then $q \in a \hjoin b$ for some $a \in X, b \in Y$. Note that $a \in X \subseteq [X\rangle$ and $b \in Y \subseteq [Y\rangle$. By the definition of type-III hyperfilter, we get $a \hjoin b \subseteq [Y\rangle$ (because $a \in \mathfrak{L}, b \in [ Y\rangle$). Similarly, we get $a \hjoin b \subseteq [X\rangle$. This gives $q \in [X\rangle \cap [Y\rangle$. Thus $B \subseteq [X\rangle \cap [Y\rangle$, and so, $[B\rangle \subseteq [X\rangle \cap [Y\rangle$.
		\item We first show that $A \subseteq \Big [[X \rangle \cup [Y\rangle \Big \rangle$. Let $a \in A$. Then $a= x \wedge y$ for some $x \in X$ and $y \in Y$. Note that $x, y \in [X \rangle \cup [Y\rangle \subseteq \Big [[X \rangle \cup [Y\rangle \Big \rangle$. So we get $a \in \Big [[X \rangle \cup [Y\rangle \Big \rangle$. Thus $A \subseteq \Big [[X \rangle \cup [Y\rangle \Big \rangle$, which implies $[A \rangle \subseteq \Big [[X \rangle \cup [Y\rangle \Big \rangle$. To show $\Big [[X \rangle \cup [Y\rangle \Big \rangle \subseteq [A \rangle$, it is enough to show that $X, Y \subseteq [A\rangle$. Let $x \in X$. Then, if $y \in Y$ is arbitrary, $a=x \wedge y \in A \subseteq [A\rangle$ and since $x \wedge y \preceq x$, we get $x \in[A\rangle$. This implies $X \subseteq [A\rangle$ which shows that $[X\rangle \subseteq [A\rangle$. By a similar argument, we get that $[Y\rangle \subseteq [A\rangle$. Therefore $[A\rangle = \Big[[X \rangle  \cup [Y\rangle \Big \rangle$. \newline For the second part, let $b \in B$. Then $b \in s \hjoin t$ for some $s \in X$ and $t \in Y$. Note that $s,t \in \langle X] \cup \langle Y] \subseteq \Big \langle \langle X] \cup \langle Y] \Big ]$. Then by the definition of hyperideal, we get $b \in s \hjoin t \subseteq \Big \langle \langle X] \cup \langle Y] \Big ]$. Thus, we get $\langle B] \subseteq \Big \langle \langle X] \cup \langle Y] \Big ]$. Now to show the reverse inclusion, we show that $X,Y \subseteq \langle B]$. Let $x \in X$. Then for any $y \in Y$, we have $x \hjoin y \subseteq B \subseteq \langle B]$. By Proposition~\ref{prop: orderproperty}, there exists $u \in x \hjoin y$ such that $x \preceq u$. Now, $u \in \langle B]$ and $x \preceq u$ implies that $x \in \langle B]$. Thus $X \subseteq \langle B]$. By a similar argument, we get $Y \subseteq \langle B]$. Hence, $\langle X] \cup \langle Y] \subseteq \langle B]$, which implies $\Big \langle \langle X] \cup \langle Y] \Big] \subseteq \langle B]$. Therefore, $\langle B]=\Big \langle \langle X] \cup \langle Y] \Big]$.
	\end{enumerate}
\end{proof}

\begin{proposition}\label{prop:6}
	Let $X$ and $Y$ be subhyperlattices of a $P$-hyperlattice $(\mathfrak{L},\hjoin\limits^P, \wedge)$, and \newline $A=X \barwedge Y=\{x \wedge y\mid x \in X, y \in Y\}, B=X \hjoin\limits^P Y=\bigcup\limits_{a \in X,b \in Y} a \hjoin\limits^P b$. Then  $\langle A] = \langle X] \cap \langle Y]$ and $[B\rangle = [X\rangle \cap [Y\rangle$.
\end{proposition}
\begin{proof}
	From Proposition~\ref{prop:4}, we have  $\langle A] \subseteq \langle X] \cap \langle Y]$ and $[B\rangle \subseteq [X\rangle \cap [Y\rangle$. Let $p \in  \langle X] \cap \langle Y]$. Then $p \preceq s$ for some $s \in h_1 \hjoin\limits^P h_2 \hjoin\limits^P \dots \hjoin\limits^P h_n$ where $h_i \in X (1 \leq i \leq n)$ and $n \in \mathbb{N}$, and $p \preceq t$ for some $t \in g_1 \hjoin\limits^P g_2 \hjoin\limits^P \dots \hjoin\limits^P g_m$ where $g_j \in Y, m \in \mathbb{N}$ and $1 \leq j \leq m$. Since $X$ and $Y$ are subhyperlattices, we get $s \in X$ and $t \in Y$. We have $p \preceq s \wedge t$ where $s \wedge t \in A \subseteq \langle A]$. This implies $p \in \langle A ]$. Thus $\langle X] \cap \langle Y] \subseteq \langle A]$. Let $b \in [X\rangle \cap [Y\rangle$. Then  $s=g_1 \wedge g_2 \wedge \dots \wedge g_k \preceq b$ for some $g_i \in X$ and $k \in \mathbb{N}$ where $1 \leq i \leq k$ and  $t=f_1 \wedge f_2 \wedge \dots \wedge f_l \preceq b$ for some $f_j \in Y$ and $l \in \mathbb{N}$ where $1 \leq j \leq l$. Since $X$ and $Y$ are subhyperlattices, we get $s \in X$ and $t \in Y$. Hence, we have $s \hjoin\limits^P t \subseteq B \subseteq [B\rangle$. Let $k_1 \in b \hjoin\limits^P b$. Then $k_1=b \vee b\vee p_1$ for some $p_1 \in P$. Note that $s \vee t \vee p_1 \preceq b \vee b \vee p_1$ where $k_2=s \vee t \vee p_1 \in s \hjoin \limits^P t$. Thus for any $k \in b \hjoin \limits^P b$, there exists $ l \in s \hjoin\limits^P t$ such that $l \preceq k$. In particular, for $b \in b \hjoin\limits^P b$, there exists $w \in s \hjoin \limits^P t$ such that $w \preceq b$. We have $w \in s \hjoin\limits^P t\subseteq [B\rangle$ and $w \preceq b$. This gives $b \in [B\rangle$. Hence $[B\rangle =[X\rangle \cap [Y\rangle$. 
\end{proof}

\begin{remark}\label{rem:counterforprop6}
	If $X$ and $Y$ are not subhyperlattices in Proposition~\ref{prop:6}, then the equality may not hold. Consider the hyperlattice given in Figure~\ref{fig10}. Let $X=\{0,b\},~ Y=\{c\}$. Clearly $X$ and $Y$ are not subhyperlattices. We have $A=X \barwedge Y = \{0\}$ and $B = X \hjoin\limits^P Y =\{c,1\}$. Note that $\langle X] =\{0,a,b \}$, $\langle Y] =\{0,a,b,c,1\}$ and $\langle A] =\{0,a\}$. Here $\langle A] =\{0,a\} \subset \{0,a,b\} = \langle X] \cap \langle Y]$.
	\begin{figure}[th]
		\begin{center}
			\resizebox{5.5cm}{!}{
				\begin{tikzpicture}
					
					\node  (h) at (0,4) {\begin{tabular}{c}
							$ \{1\} $\\$1 $
					\end{tabular}};
					
					\node  (b) at (-3,2) {\begin{tabular}{c}
							$ \{b\} $\\$b $
					\end{tabular}};;

					\node  (a) at (-3,0) {\begin{tabular}{c}
							$ \{a\} $\\$a $
					\end{tabular}};
					
					\node  (c) at (3,0.8) {\begin{tabular}{c}
							$ \{c,1\} $\\$c$
					\end{tabular}};
					\node  (0) at (0,-2) {\begin{tabular}{c}
							$ \{0,a\} $\\$0 $
					\end{tabular}};
					\draw[thick]  (h)--(b)--(a)--(0)--(c)--(h);
					
				\end{tikzpicture}
			}
			\caption{Hyperlattice considered in Remark~\ref{rem:counterforprop6}}\label{fig10}
		\end{center}
	\end{figure}
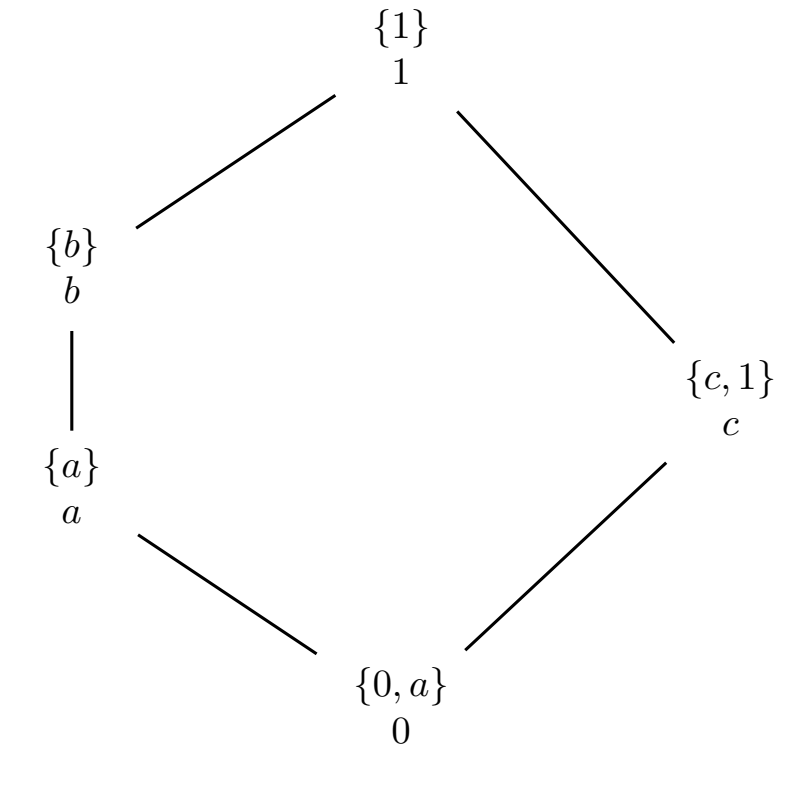
\end{remark}

The expressions for meet and join in $(\CS(\mathfrak{L}),\leq_b)$ given by Lavanya and Parameshwara Bhatta~\cite{lavanya1996new} motivated Bhargava to develop analogous results in trellises \cite{Bhargava2026}. The result below on join hyperlattices is derived from the ideas in these works.
\begin{theorem}
	Let $X$ and $Y$ be any two convex subhyperlattices of a $P$-hyperlattice $\mathfrak{L}$, and $A=X \barwedge Y=\{x \wedge y\mid x \in X, y \in Y\}, B=X \hjoin\limits^P Y=\mathlarger{\bigcup\limits_{a \in X,b \in Y} a \hjoin\limits^P b}$. Then $X \wedge Y =\langle A \rangle$ and $X \vee Y =\langle B \rangle$ where $X \wedge Y$ and $X \vee Y$ are respectively the meet and join of $X$ and $Y$ in $(\CS(\mathfrak{L}),\sqsubseteq)$.
\end{theorem}
\begin{proof}
	This follows from Proposition~\ref{idealorderlattice},  Proposition~\ref{prop:4} and Proposition~\ref{prop:6}.
\end{proof}

\begin{definition}[\cite{Gratzer}]
	A subposet $Q$ of a poset $P$ is cover preserving if whenever $x \prec y$ in $Q$ for any $x,y \in Q$, then $x \prec y$ in $P$.
\end{definition}

\begin{theorem}\label{thm:sublatticeK}
	For a $P$-hyperlattice $\mathfrak{L}$, the map $f: \CS(\mathfrak{L}) \rightarrow I(\mathfrak{L}) \times F(\mathfrak{L})$ defined by, for $C \in \CS(\mathfrak{L})$, $f(C)=(\langle C] , [C \rangle)$, is an injective homomorphism. Moreover, the sublattice $\mathcal{K}=\{ (I,D): I \in I(\mathfrak{L}), D\in F(\mathfrak{L}), I \cap D \neq \emptyset\}$ of $I(\mathfrak{L}) \times F(\mathfrak{L})$  is cover preserving and is  isomorphic to $(\CS(\mathfrak{L}),\sqsubseteq)$. 
\end{theorem}
\begin{proof}
	The map $f$ is well defined. Also, $f$ is a one-to-one map. Note that for any $X_1, Y_1 \in \CS(\mathfrak{L})$, by Remark~\ref{rem:idealbyCSH}, we have $\langle X_1 \wedge Y_1] =\langle X_1 ] \wedge \langle Y_1]$ and $[X_1 \wedge Y_1 \rangle =\Big [[X_1 \rangle \cup [Y_1\rangle \Big \rangle =[X_1 \rangle \wedge [Y_1\rangle$. Now,
	\begin{eqnarray*}
		f(X_1 \wedge Y_1) &=&\big (\langle X_1 \wedge Y_1],[X_1 \wedge Y_1 \rangle \big) \\
		&=& \big (\langle X_1] \wedge \langle Y_1], [X_1\rangle \wedge [Y_1\rangle \big)\\
		&=& \big(\langle X_1], [X_1\rangle\big) \wedge \big(\langle Y_1], [Y_1\rangle\big)\\
		&=& f(X_1) \wedge f(Y_1)
	\end{eqnarray*} 
	Therefore, $f$ is a meet homomorphism of lattices. Similarly, $f$ is a join homomorphism of lattices. Thus, $f$ is a homomorphism. Clearly, $f(\CS(\mathfrak{L}))=\mathcal{K}$ and hence $\CS(\mathfrak{L}) \cong \mathcal{K}$. Now, we show that $\mathcal{K}$ is cover preserving. Let $(I_1,D_1), (I_2,D_2) \in \mathcal{K}$ such that $(I_1 ,D_1) \prec (I_2,D_2)$ in $\mathcal{K}$. It is enough to show that $I_1=I_2$ and $D_1 \prec D_2$ in $F(\mathfrak{L})$ or $I_1 \prec I_2$ in $I(\mathfrak{L})$ and $D_1=D_2$. Suppose $I_1 \neq I_2$ and $D_1 \neq D_2$. Then $I_1 \subsetneq I_2$ and $D_1 \supsetneq D_2$. Let $p \in I_2 \setminus I_1$. Take $I=I_1 \vee \langle p]$ and $D=D_1$. Then $(I,D) \in \mathcal{K}$ and $(I_1,D_1) < (I,D) < (I_2,D_2)$ in $\mathcal{K}$. This is a contradiction to the fact that $(I_2, D_2)$ covers $(I_1,D_1)$ in $\mathcal{K}$. Thus $I_1=I_2$ or $D_1=D_2$. Suppose $I_1=I_2$. If $D_1 < D < D_2$ in $F(\mathfrak{L})$ for some $D\in F(\mathfrak{L})$, then, $(I_2,D) \in \mathcal{K}$ and $(I_1,D_1) < (I_2,D) < (I_2,D_2)$ in $\mathcal{K}$, which is a contradiction. Hence, $D_1 \prec D_2$ in $F(\mathfrak{L})$. Similarly, if $D_1=D_2$, then $I_1 \prec I_2$ in $I(\mathfrak{L})$.
\end{proof}

\begin{remark}
	The sublattice $\mathcal{K}$ of $I(\mathfrak{L})\times F(\mathfrak{L})$ in Theorem~\ref{thm:sublatticeK} may not be convex. For example, consider the hyperlattice $\mathfrak{L}$ given in Example~\ref{eg1}. We note that $(\{0,c\},\{c,f,h,1\}) \leq (\{0,c\}, \{h,1\}) \leq (\{0,c,d,e,h\},\{h,1\})$ and  $(\{0,c\}, \{h,1\}) \in I(\mathfrak{L}) \times F(\mathfrak{L})$ but $(\{0,c\}, \{h,1\}) \notin \mathcal{K}$.
\end{remark}
\begin{proposition}\label{prop:7}
	For a $P$-hyperlattice $\mathfrak{L}$, $I(\mathfrak{L})$ is the principal ideal of $(\CS(\mathfrak{L}),\sqsubseteq)$ generated by $\mathfrak{L}$ and $F(\mathfrak{L})$ is the principal filter  of $(\CS(\mathfrak{L}),\sqsubseteq)$ generated by $\mathfrak{L}$ .
\end{proposition}
\begin{proof}
	Let $X \in \CS(\mathfrak{L})$ such that $X \sqsubseteq \mathfrak{L}$. Then $[X\rangle \supseteq [\mathfrak{L}\rangle =\mathfrak{L}$. Thus $[X\rangle =\mathfrak{L}$. This shows that $X=\langle X] \cap [X\rangle =\langle X] \cap \mathfrak{L} =\langle X] \in I(\mathfrak{L})$. Suppose $X \in I(\mathfrak{L})$, then $X = \langle X] \cap \mathfrak{L}$. Being a hyperideal, $X$ is a convex subhyperlattice. By the uniqueness of the representation for $X$, we get $[X\rangle =\mathfrak{L}$. Hence we get that $X \sqsubseteq \mathfrak{L}$. By a similar argument, we can show that $F(\mathfrak{L})$ is the principal filter generated by $\mathfrak{L}$. 
\end{proof}

\begin{definition}[\cite{Gratzer2011}]
	Let $p_i=q_i$ be identities for $i \in I$. The class $\kappa$ of all lattices satisfying the identities $p_i=q_i$, for $i \in I$, is called a equational class of lattices.
\end{definition}
\begin{corollary}\label{cor:equationalclass}
	Let $\mathfrak{L}$ be a $P$-hyperlattice and $\gamma$ be an equational class of lattices. If $I(\mathfrak{L})$ and $F(\mathfrak{L})$ lie in $\gamma$, then $(\CS(\mathfrak{L}),\sqsubseteq)$ lies in $\gamma$. Conversely, if $(\CS(\mathfrak{L}),\sqsubseteq)$ lies in $\gamma$, then both $I(\mathfrak{L})$ and $F(\mathfrak{L})$ lie in $\gamma$.
\end{corollary}
\begin{proof}
	This follows from Theorem~\ref{thm:sublatticeK} and Proposition~\ref{prop:7}.
\end{proof}
In general, $I(\mathfrak{L})$ and $F(\mathfrak{L})$ need not be in a same equational class. For instance, in Example~\ref{eg1}, $I(\mathfrak{L})$ is modular, whereas $F(\mathfrak{L})$ is non-modular. Secondly, even when $I(\mathfrak{L})$ and $F(\mathfrak{L})$ are in the same equational class,  $(\CS(\mathfrak{L}) \cup \{\emptyset\}, \subseteq)$ may fail to lie in the same class. In the $P$-hyperlattice shown in Figure~\ref{fig13}, $I(\mathfrak{L}), F(\mathfrak{L})$ and $(\CS(\mathfrak{L}), \sqsubseteq)$ are distributive whereas $(\CS(\mathfrak{L}) \cup \{\emptyset\}, \subseteq)$ is not distributive. Moreover, the result stated in Corollary~\ref{cor:equationalclass} does not hold if we consider  $(\CS(\mathfrak{L}) \cup \{\emptyset\}, \subseteq)$ instead of $(\CS(\mathfrak{L}), \sqsubseteq)$. 
\begin{example}\label{eg9}
	Let $\mathfrak{L}=\{0,p,q,r,s,1\}$. The $P$-hyperlattice with $P=\{0,q\}$ is given in Figure~\ref{fig13}.
	\begin{figure}[th]
		\begin{center}
			\resizebox{5.5cm}{!}{
				\begin{tikzpicture}
					
					\node  (1) at (0,4) {\begin{tabular}{c}
							$ \{1\} $\\$1 $
					\end{tabular}};
					
					\node  (r) at (-3,2) {\begin{tabular}{c}
							$ \{r,1\} $\\$r $
					\end{tabular}};;

					\node  (p) at (-3,0) {\begin{tabular}{c}
							$ \{p,s\} $\\$p $
					\end{tabular}};
					
					\node  (q) at (3,0) {\begin{tabular}{c}
							$ \{q\} $\\$q $
					\end{tabular}};
					
					\node  (s) at (3,2) {\begin{tabular}{c}
							$ \{s\} $\\$s$
					\end{tabular}};
					\node  (0) at (0,-2) {\begin{tabular}{c}
							$ \{0,q\} $\\$0 $
					\end{tabular}};
					\draw[thick]  (1)--(r)--(p)--(0)--(q)--(s)--(1);
					\draw[thick]  (p)--(s);
					
				\end{tikzpicture}
			}
			\caption{$P$-hyperlattice considered in Example~\ref{eg9}}\label{fig13}
		\end{center}
	\end{figure}
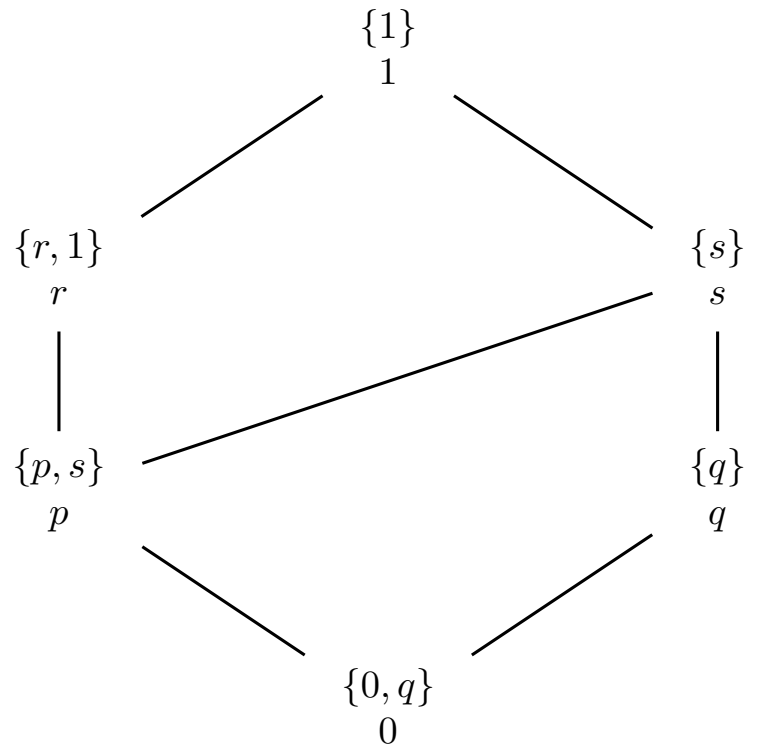
\end{example}

Let $\mathcal{P}^*(\mathfrak{L})$ denote the collection of all non-empty subsets of a $P$-hyperlattice $\mathfrak{L}$. For $P_1, P_2 \in \mathcal{P}^*(\mathfrak{L})$, define $P_1 \leq P_2$ if and only if $\langle P_1] \subseteq \langle P_2]$ and $[P_1\rangle \supseteq [P_2\rangle$. Then clearly $\leq$ is reflexive and transitive. Now, define a new relation $\theta$ on $\mathcal{P}^*(\mathfrak{L})$ as follows: For $P_1, P_2 \in \mathcal{P}^*(\mathfrak{L})$, define $P_1 \equiv P_2 (\theta)$ if and only if $P_1 \leq P_2$ and $P_2 \leq P_1$. We note that $\theta$ is an equivalence relation on $\mathcal{P}^*(\mathfrak{L})$. Note that $P_1 \equiv \langle P_1 \rangle (\theta)$, and $P_1 \equiv P_2 (\theta)$ if and only if $\langle P_1 \rangle =\langle P_2 \rangle$. Hence the relation $\leq$ defined on $\mathcal{P}^*(\mathfrak{L}) / \theta$ by, for $[P_1], [P_2] \in \mathcal{P}^*(\mathfrak{L}) / \theta,$ $[P_1]\leq [P_2]$ if and only if $\langle P_1 \rangle \leq \langle P_2 \rangle$ is a partial order on $\mathcal{P}^*(\mathfrak{L}) / \theta$. The map $[P_1] \rightarrow \langle P_1 \rangle$ defines an order isomorphism between the two posets $(\mathcal{P}^*(\mathfrak{L}) / \theta, \leq)$ and $(\CS(\mathfrak{L}), \sqsubseteq)$. Since $(\CS(\mathfrak{L}), \sqsubseteq)$ is a lattice, we note that $(\mathcal{P}^*(\mathfrak{L}) / \theta, \leq)$ is also a lattice. \newline
Thus we have the following proposition.
\begin{proposition}
	For a $P$-hyperlattice $\mathfrak{L}$,
	\begin{equation*}
		(\CS(\mathfrak{L}), \sqsubseteq) \cong (\mathcal{P}^*(\mathfrak{L}) / \theta, \leq)
	\end{equation*}
\end{proposition}

The notion of convexity is well understood in lattices and extends naturally to some classes of strong hyperlattices, but it fails to exhibit the same structural properties in an arbitrary join hyperlattice. In particular, a convex subhyperlattice need not admit a unique representation in terms of hyperideals and hyperfilters. In addition to this, several results valid in $P$-hyperlattices and Nakano hyperlattices fail to extend directly. To overcome this, we introduce the concept of strongly convex sets, which refines the existing notion of convexity while preserving its properties. This concept serves as a natural generalization, ensuring that results established for $P$-hyperlattices remain valid in the broader setting of arbitrary hyperlattices.
\section{Strongly convex sets in $\mathfrak{L}$}
In this section, we introduce the notion of strongly convex sets in an arbitrary join hyperlattice and study their basic properties. We show that every strongly convex set $X$ admits a unique representation in terms of $\langle X]$ and $[X\rangle$. This characterization allows us to extend several structural results known for $P$-hyperlattices to a more general setting. 
\begin{definition}
	A non-empty set $S$ of $\mathfrak{L}$ is said to be strongly convex if $S=I \cap F$ for some hyperideal $I$ and type-III hyperfilter $F$.
\end{definition}

\begin{remark}
	By Theorem~\ref{thmC}, every strongly convex subset of $\mathfrak{L}$ is a convex subhyperlattice, but the converse does not hold in general. This is illustrated in Example~\ref{eg10}, where $\{x_1\}$ is a convex subhyperlattice which is not strongly convex. The example also illustrates that $\{1\}$ is a convex subhyperlattice which is a strongly convex set.  Furthermore, in a $P$-hyperlattice or a Nakano hyperlattice, every convex subhyperlattice is strongly convex.
\end{remark}
We denote the set of all strongly convex sets of $\mathfrak{L}$ by $\SC(\mathfrak{L})$.
\begin{proposition}
	Let $X$ be a strongly convex set of $(\mathfrak{L}, \hjoin, \wedge)$. Then $X=\langle X] \cap [X\rangle$.
\end{proposition}
\begin{proof}
	Let $X=I \cap F$ for some  hyperideal $I$ and type-III hyperfilter $F$. Then, clearly, $\langle X] \subseteq I$ and $[X\rangle  \subseteq F$. Let $I_1$ be a hyperideal of $\mathfrak{L}$ containing $X$. Suppose that $a \in I$ and let $x \in X$ be arbitrary. Then, $a \hjoin x \subseteq I$. Also, $x \in F, a \in \mathfrak{L}$ implies that $a \hjoin x \subseteq F$ since $F$ is a type-III hyperfilter. This implies $a \hjoin x \subseteq I \cap F = X \subseteq I_1$. Note that there exists $u \in a \hjoin x$ such that $a \preceq u$. Thus, $u \in I_1$ and $a \preceq u$ implies that $a \in I_1$. Hence, $I \subseteq I_1$. Now $\langle X] $ is a hyperideal containing $X$. By the above argument, we get that $I \subseteq \langle X]$, and hence $I=\langle X]$. Now, let $F_1$ be any type-III hyperfilter of $\mathfrak{L}$ containing $X$. Let $p \in F$. Then, for any $q \in X$, we have $p \wedge q \in F$. Also, since $q \in I$, we get that $ p \wedge q \in I$. Therefore $p \wedge q \in X \subseteq F_1$. Now, $p \wedge q \in F_1$ and $p \wedge q \preceq p$ implies that $p \in F_1$. This shows that $F \subseteq F_1$. Therefore $F=[X\rangle$.
\end{proof}
The following theorems are generalizations of Proposition~\ref{proposition:idealorder}, Proposition~\ref{idealorderlattice} and  Theorem~\ref{thm:sublatticeK} to arbitrary join hyperlattices. Their proofs are omitted, as they follow directly from the arguments used in the corresponding results. In the sequel, let $\mathfrak{L}$ be such that $\SC(\mathfrak{L})=\CS(\mathfrak{L})$.

\begin{theorem}
	For a join hyperlattice $\mathfrak{L}$, $(\SC(\mathfrak{L}),\sqsubseteq)$ is a poset.
\end{theorem}
\begin{theorem}
	The poset $(\SC(\mathfrak{L}),\sqsubseteq)$ is a lattice with meet and join as given in Proposition~\ref{idealorderlattice}.
\end{theorem}
\begin{theorem}
	For a hyperlattice $\mathfrak{L}$, the map $f: \SC(\mathfrak{L}) \rightarrow I(\mathfrak{L}) \times F(\mathfrak{L})$ defined by, for $C \in \SC(\mathfrak{L})$, $f(C)=(\langle C] , [C \rangle)$, is an injective homomorphism. Moreover, the sublattice $\mathcal{K}=\{ (I,D): I \in I(\mathfrak{L}), D\in F(\mathfrak{L}), I \cap D \neq \emptyset\}$ of $I(\mathfrak{L}) \times F(\mathfrak{L})$  is cover preserving and is  isomorphic to $(\SC(\mathfrak{L}),\sqsubseteq)$ . 
\end{theorem}

\section{Conclusion}

It is well known that the family $\CS(\mathfrak{L}) \cup \{\emptyset\}$, ordered by inclusion, forms a lattice. In this work, we have extended the study of this structure. First, we have shown that $\CS(\mathfrak{L}) \cup \{\emptyset\}$ is not only a lattice but an algebraic lattice, and we have characterized its compact elements. Then, we have introduced two new orders on $\CS(\mathfrak{L})$, namely $\leq_b$ and $\sqsubseteq$, and have examined their properties. We have proved that for $A, B \in \CS(\mathfrak{L})$, $A \leq_b B$ implies $A \sqsubseteq B$. However, these two orders are not generally equivalent, though they coincide for $P$-hyperlattices. We have shown, by example, that in an arbitrary join hyperlattice, the relation $\sqsubseteq$ need not be a partial order. Next, we have characterized convex subhyperlattices of $P$-hyperlattices using initial and final segments, and have embedded the lattice $(\CS(\mathfrak{L}), \sqsubseteq)$ into the product lattice $I(\mathfrak{L}) \times F(\mathfrak{L})$. We have introduced the notion of strongly convex sets and have shown that every such set is a convex subhyperlattice, but not vice versa. Although many results are given for $P$-hyperlattices, the key arguments do not rely on the $P$-hyperlattice condition. Consequently, these results also hold for hyperlattices satisfying $\SC(\mathfrak{L})=\CS(\mathfrak{L})$. Finally, we have justified the use of type-III hyperfilters in defining strongly convex sets, since type-III hyperfilters, unlike other types, form subhyperlattices. We have proved that when $\mathfrak{L}$ is a $P$-hyperlattice, $(\CS(\mathfrak{L}), \leq_b)$ is a lattice, and we have given an example of a non-strong hyperlattice where this does not hold. Based on our observations, we have conjectured that $(\CS(\mathfrak{L}), \leq_b)$ is a lattice for every strong hyperlattice. In conclusion, it is natural to investigate whether the operator $\SC$ preserves finite products of hyperlattices and to identify the structural conditions under which such preservation holds. Additionally, examining how homomorphisms and congruence relations on a hyperlattice $\mathfrak{L}$ carry over to the induced structure on $\SC(\mathfrak{L})$ offers scope for further research.

\section*{Acknowledgement(s)}

The first author acknowledges the Department of Science and Technology (DST), Government of India, for the DST INSPIRE Fellowship to carry-out this research work. The first, fourth, and sixth authors acknowledge the Manipal Institute of Technology, Manipal Academy of Higher Education, Manipal for its kind encouragement. The second author acknowledges the Department of Mathematics, Mangalore university, Mangalagangotri for its support. The third and fifth authors thank the Manipal Institute of Technology Bengaluru, Manipal Academy of Higher Education, Manipal for its encouragement and institutional support.

\bibliographystyle{plain}
\bibliography{bibreport}

\textit{Preprint of an article submitted for consideration in Journal of Algebra and its Applications $\copyright$ $2026$ World Scientific Publishing Co Pte Ltd \url{https://www.worldscientific.com/worldscinet/jaa}.}
\end{document}